\documentclass{amsart}
\usepackage{color}
\usepackage[toc,page]{appendix}
\usepackage[utf8]{inputenc}
\usepackage[margin=1in,left=1in]{geometry}
\usepackage{setspace}
\usepackage{graphicx}
\usepackage{amssymb, amsmath, amsthm, graphics}
\usepackage{mathabx,epsfig}
\usepackage[mathscr]{euscript}
\usepackage{tikz}
\usetikzlibrary{calc}
\usetikzlibrary{matrix}
\usepackage{hyperref}     

\usepackage{latexsym}
\usepackage[all]{xy}
\usepackage{color}

\usepackage{bbold}

\usepackage{tikz}
\input xy
\xyoption{all}
\usepackage{euscript}
\usepackage{multirow}
\usetikzlibrary{positioning}
\usetikzlibrary{shapes.geometric}
\usetikzlibrary{shapes.misc}
\usetikzlibrary{calc}
\usetikzlibrary{positioning}

\usepackage{pgflibraryarrows}
\usepackage{pgflibrarysnakes}

\usetikzlibrary{trees} 
\usetikzlibrary[trees] 

\newtheorem{theorem}{Theorem}
\newtheorem{definition}[theorem]{Definition}
\newtheorem{lemma}[theorem]{Lemma}
\newtheorem{corollary}[theorem]{Corollary}
\newtheorem{proposition}[theorem]{Proposition}

\newtheorem{remark}{Remark}

\newtheorem{example}{Example}

\numberwithin{equation}{section}

\renewcommand{\(}{\begin{equation*}}
\renewcommand{\)}{\end{equation*}}
\newcommand{\bea}{\begin{eqnarray*}}
\newcommand{\eea}{\end{eqnarray*}}

\def\proof {{Proof.}\hspace{7pt}}

\def\endofproof {\hfill{$\Box$}\\}

\newcommand{\Or}{{\rm O}}
\newcommand{\Spin}{{\rm Spin}}
\newcommand{\String}{{\rm String}}
\newcommand{\Five}{{\rm Fivebrane}}
\newcommand{\Nine}{{\rm Ninebrane}}
\newcommand{\psix}{\tfrac{1}{6}p_2}

\newcommand{\beq}{\begin{equation}}
\newcommand{\eeq}{\end{equation}}

\newcommand{\op}[1]{\ensuremath{\operatorname{#1}}}

\newcommand{\theproof}{\noindent {\bf Proof.\ }}

\numberwithin{equation}{section}

\renewcommand{\(}{\begin{equation}}
\renewcommand{\)}{\end{equation}}

\newcommand{\Exterior}{\scalebox{.8}{\ensuremath \bigwedge}}

\def\1{{\bf 1}}

\def\<{\langle}
\def\>{\rangle}

\numberwithin{equation}{section}


\newcommand{\R}{\ensuremath{\mathbb R}}

\newcommand{\ZZ}{\ensuremath{\mathbb Z}}
\newcommand{\Z}{\ensuremath{\mathbb Z}} 
\newcommand{\Q}{\ensuremath{\mathbb Q}}

\newcommand{\CC}{\ensuremath{\mathbb C}}

\author{Hisham Sati and Matthew Wheeler}
\email{hsati@nyu.edu, mgwheeler@math.arizona.edu}
\title{Variations of rational higher tangential structures}

\begin{document}
\maketitle 

\begin{abstract}
The study of higher tangential structures, arising from higher connected
covers of Lie groups (String, Fivebrane, Ninebrane structures),
require considerable machinery for a full description, especially 
for connections to geometry and applications.
With utility in mind, in this paper we study these structures at the rational
level and by considering Lie groups as a starting point for defining each 
of the higher structures, making close connection to $p_i$-structures. 
We indicatively call these (rational) {\it Spin-Fivebrane} and 
{\it Spin-Ninebrane} structures. 
We study the space of such structures and  characterize their variations, 
which reveal interesting effects whereby variations 
of higher structures are arranged to systematically involve lower ones. 
We also study the homotopy type of  the gauge 
group corresponding to bundles equipped with the higher rational structures
that we define. 

\end{abstract}

\tableofcontents

\section{Introduction}
Manifolds have been classically studied through structures associated with their 
tangent bundles leading to characterizations via obstruction theory and characteristic classes
\cite{Steenrod}\cite{MS}\cite{Hu}. 
We examine the rationalization of tangential structures, with 
an emphasis on structures arising from higher connected covers of Lie groups.  
That is, we consider rationalizing the higher structure groups and their classifying spaces,
namely String \cite{Kil}\cite{ST}, Fivebrane \cite{SSS2}\cite{SSS3}, and Ninebrane 
structures \cite{9brane}. This has a simplifying feature in that tangential structures 
from obstruction theory \cite{Steenrod}\cite{Hu} are algebraically placed 
in the setting of rational homotopy theory  \cite{FHT}\cite{FHT2}\cite{FOT}\cite{GM}\cite{BG}.   
This setting allows us to filter out the torsion in our spaces
thereby enabling us to have a much better handle on 
 some aspects of these otherwise 
formidable structures. However, on the flip side, a complication arises when 
wishing to describe the rationalizations as spaces, since localization in general 
give rise to topological spaces which are not always nice
\cite{Dr}\cite{Ne}\cite{HMR}\cite{BK}. 
Our discussion will strike a balance between the two competing aspects and
our goal in this paper is to highlight those features that have transparent descriptions. 

\medskip
Another aim of this paper is to investigate to which extent one can make use of 
the more familiar Lie group structures 
in describing the higher ones. In the standard Whitehead tower construction \cite{Wh}, 
a structure at a given level is built from
 the structure at the preceding level. However,
as we go up in levels, the difficulty in describing the structures in an explicit manner
which is  amenable to (higher) geometry and to applications seems to grow considerably. 
Therefore, it  would be desirable to explore how much of  the bundles 
of higher connected covers can be described 
using the Spin group (being a Lie group) rather than having, for instance,
 to go through String to describe Fivebrane and through Fivebrane to describe 
 Ninebrane structures and so on. Of course one can deal with these structures directly
 (see \cite{SSS3}\cite{FSSt}\cite{FSS1}\cite{FSS2}\cite{9brane}), but 
 that requires considerable machinery. 
 Here we instead take a step back and aim to explore to which extent more classical 
 techniques can be used to probe these structures.

 \medskip 
It turns out that the two features, namely rationalizing and utilization of Lie groups in describing the higher connected covers, go hand in hand. 
The purpose of this paper is to provide a straightforward such description.  
One of the useful results which makes this possible is that of 
Neisendorfer \cite{Ne} (see also \cite{MM}) which states that every finite 2-connected complex
can be rationally recovered from its $n$-connected cover 
for any $n$. This is much stronger than saying that these spaces must 
have nontrivial homotopy in infinite dimensions. As explained in \cite{Dr} 
it says that this `infinite tail' has all the information needed to reconstruct the 
`lower-dimensional information'. This then allows us to appropriately introduce 
(in Sec. \ref{Sec var5}  and Sec. \ref{Sec var9}) 
the notions of rational {\it Spin-Fivebrane} and {\it Spin-Ninebrane structures} as the
desired  structures
arising from starting with Spin rather than with String and Fivebrane,
respectively.

\medskip
In general, we would like to start with a Lie group $G$ and then rationalize, via the rationalization 
or localization at $\Q$ functor $L_\Q$, as well as take connected covers 
simultaneously. A natural question then is 
whether these operations are compatible, in the sense of the existence of a diagram
of the form 
$$
\xymatrix@R=1.8em{
X \ar[rr]^{\langle n \rangle} \ar[d]_{L_\Q}
&& X \langle n \rangle \ar[d]^{L_\Q}
\\
X_\Q \ar[rr]^{?} && ~X_\Q\langle n \rangle\;.
}
$$
One of the simplifying features of the
 process of rationalization is that for group-like spaces  it
has the effect of killing off any nonabelian structure that exists. Schematically,
\begin{center}
Rationalization $\leadsto$ Homotopy abelianization.
\end{center}
This then gives that all connected cover groups will not only have rational models,
but that these will be homotopy abelian. Corresponding statements about 
the classifying spaces are deduced similarly. 
Our main focus will be on the secondary structures arising from the groups in the Whitehead 
tower of the orthogonal group $O(n)$ and, in particular, we will focus on the rationalizations of these groups.  
Given an ${\rm O}(n)\langle k\rangle$-bundle $P\to M$, the obstruction to lifting  the structure 
group to the $k$-connected cover ${\rm O}(n)\langle k+1\rangle$ is given by a cohomology 
class on $M$ obtained by pulling back the generator 
$\theta_{k+1}\in H^{k+1}(B({\rm O}(n)\langle k\rangle);\pi_k({\rm O}(n)))$ along the classifying 
map $f:M\to B({\rm O}(n)\langle k\rangle)$.   

%

\medskip
Note that for Lie groups, maps between their classifying spaces
can be determined via Lie theory, with an intimate connection to rational cohomology 
\cite{AM}. In fact, homomorphisms $H^*(BG; \Q) \to H^*(BG';\Q)$ 
 determine corresponding homomorphism with coefficients in $\Z_{(p)}$,
 the ring of integers localized at a prime $p$, and in $\Z/p$ the field 
 of integers modulo $p$, except for a finite number of primes \cite{AM}. This 
 indicates that rational cohomology knows quite a bit about the structure of the classifying spaces.
 We hope that our investigation on rational cohomology of classifying spaces of the 
 connected covers will eventually carry some of the similar features.

 \medskip
 From a third point of view, we are interested in considering 
 classes in degrees 3, 7, and 11 arising from the fibers of the bundles, 
 rather than directly from classifying spaces. These secondary classes 
 emerge when trivializing the topological obstructions that occur in degrees
  4, 8 and 12. We are interested in the variation of the structures, i.e. 
  in the space of such structures. Given that we are considering  higher structures
 in  a way which builds on  all the lower levels that precede it, 
 we ask how the variation of that top structure depends on the 
 lower ones all the way down to  the bottom-most level, which is a Spin structure.

\medskip
The rational structures we consider will be 
characterized by trivializations of (fractions of) the rational  Pontrjagin
 classes $p_i$. Indeed, we can alternatively view the above 
 structures as variants of $p_i$-structures, so that we have lifting diagrams 
$$
{\small
\xymatrix{
&& BO \langle p_i \rangle \ar[d] \\
X \ar@{..>}[urr]^-{\overline{f}} \ar[rr]^-{f}  && ~BSO=BO\langle w_1 \rangle\;,
}
}
$$
where $w_1$ is the first Stiefel-Whitney class.
Note that we have the obstructions as $\tfrac{1}{m}p_i^\Q$, whose vanishing 
is equivalent to the vanishing of the rational Pontrjagin class $p_i^\Q$ itself,
due to the absence of any torsion.
The notion of  $p_1$-structures plays an important role in quantum field theories on
extended surfaces and 3-manifolds  \cite{BHMV}\cite{Seg}\cite{BN}\cite{FSV}.
 Extension to Fivebrane and Ninebrane structures has been considered 
 in \cite{9brane}, and the corresponding twisted versions come up in \cite{top}. 
We will be interested also in the finite rank case as well as the indefinite signature case,
where our main objects will be covers of $SO(n)$ and $SO(q, n-q)$ discussed in Sec. 
\ref{Sec rank}.

\medskip
The rational cohomology of $BSO(n)$ splits into cases according to 
whether $n$ is even or odd
$$
H^*(BSO(2n+1); \Q)=\Q[p_1, \cdots, p_n]\;,
\qquad
H^*(BSO(2n); \Q)=\Q[p_1, \cdots, p_n, e]/(e^2-p_n)\;,
$$
where $e \in H^{2n}(BSO(2n); \Q)$ is the Euler class
and $p_i\in H^{4i}(BSO(2n); \Q)$ are the rational 
Pontrjagin classes. 
Unlike the integral versions, the rational Pontrjagin classes 
are topological invariants \cite{Nov}\cite{RW}, which makes them 
reliable under homeomorphisms. A vanishing criterion for  $p_1^\Q$ 
 is given in \cite{Ho}.
 It is a classic result (see \cite{We}) that restrictions of the classes
$p_i$ to the classifying space of finite-dimensional vector bundles
satisfy the vanishing relations
$
p_{n+k}=0 \in H^{4n+4k}(BO(2n); \Q) 
$
 for $k>0$.

\medskip
 The rational Pontrjagin classes 
 have been 
 used in \cite{EOW} in the context of
cobordism spectra. 
A version of the Witten genus can be described by requiring the rational 
first Pontrjagin class to vanish; see e.g. \cite{De}\cite{CHZ} (and references 
therein), where a similar definition of a rational structure is used. There, a Spin 
manifold $M$ is a rational $BO \langle 8 \rangle$ manifold if and only if $p_1(M)$ is a torsion
class.  Furthermore, the rational Pontrjagin classes are used in classifying bundles 
 in \cite{Kr1}, where it is shown that  rank $4n$
 vector bundles over the $4n$-sphere $S^{4n}$ are classified by their 
 Euler class and the rational $n$th Pontrjagin class $p_n^\Q$ for 
 $n=1,2$. Such bundles classically arise in determining obstructions to lifting to 
 higher connected covers (see \cite{SSS2}). 
We, therefore, consider the question of the relation of 
the connectivity to the rank in Section \ref{Sec rank}. 

\medskip
The rational cohomology of the String group has been considered 
in \cite{SSS2}\cite{BS}. 
Also in specific ranks in relation to connectivity degree,
 $BO (2n) \langle n \rangle$ appear in the context of cobordism categories 
 \cite{EOW}, 
 where the isomorphism 
 $
 H^*(BO(2n) \langle n \rangle; \Q)[-2n]\cong H^*(MT\theta^n; \Q)
 $
 of graded vector spaces 
 is established. Here $MT\theta^n$ is the Madsen-Tillman cobordism spectrum 
with a tangential structure $\theta^n$, i.e. a structure on a space
associated to an $n$-connected cover. For $n=4, 8$ and 12, this 
corresponds in our terminology 
 to rational String, Fivebrane, and Ninebrane structures, 
respectively.

\medskip
We consider minimal models  (see \cite{FHT}\cite{FHT2}\cite{FOT}\cite{GM}\cite{BG})  
 for our rational connected covers straightforwardly  in Sec. \ref{Sec min}.
The main idea of Sullivan's approach to rational homotopy theory is to create a functor from the category of 1-connected topological spaces to the category of differential graded commutative  algebras (DGCAs) over $\Q$.
%
%
Such a DGCA is of the form $(\Exterior V, d)$ where the underlying algebra is free commutative and such that there is a basis which admits an ordering so that $d(x_\alpha)\in\Exterior(x_\beta)_{\beta<\alpha}$. Furthermore, $(A, d)$  is {minimal} if the image of the differential $d$ is contained in the set of decomposable elements, and a {minimal model} is a quasi-isomorphism $\varphi:(\Exterior V,d)\to(A,d)$ where $(\Exterior V,d)$ is a minimal Sullivan DGCA.  In fact every DGCA has a minimal model, and this model is unique up to isomorphism.  

\medskip
Then in the following sections we consider variations 
on rational String and Fivebrane structures. We use the word variation to 
mean two things at the same time: First, that we consider variations on the 
notion of Fivebrane and Ninebrane structures. Second, we consider variations 
of the actual structures (in their `parameter space') and consider how these are given in 
terms of structures stemming from lower levels in the Whitehead tower. 
 Theorem \ref{rfiv} demonstrates the degree to which the underlying Spin bundle 
 can be used to classify lifts of the String bundles rationally. By defining these classes 
 via their restriction on each fiber, we have many nice parallels between the integral and rational cases, as well as between those classes defined on the Spin bundle and those on the String bundle.  In some sense, we find that rationally all the information for Fivebrane and Ninebrane structures is encoded in the underlying Spin bundles.  Similar arguments and results hold for the Spin-Ninebrane case, except now variations of these involve both String and Fivebrane 
  structure classes. 
  
\medskip  
  In Sec. \ref{Sec gauge}, we consider automorphisms of the rational structures 
that we introduced. Rational automorphisms of fiber bundles are considered generally
in \cite{Smith}.  Since our higher groups are rationally abelianized, the description 
will be more straightforward, making use of classic results on mapping spaces 
to Eilenberg-MacLane spaces  \cite{Haef}\cite{Hansen}\cite{Thom}. We also 
study the connected covers of the gauge group $\mathcal{G}$  itself and study when 
a variant of String, Fivebrane, and Ninebrane lifts of $\mathcal{G}$ are possible in relation 
to corresponding lifts of the structure group $G$. This turns out to impose strong constraints 
on both on the underlying space $X$ as well as on $G$ in  a correlated manner. 
We make use of the results on
rational homotopy of mapping spaces in \cite{FO}, which generalize
those of \cite{Wo}.

\medskip
Some of the calculations in this note are based on the second author's PhD 
thesis \cite{Matt}. The point of view and constructions  developed here  naturally 
lead to connections to geometry, which we will leave for a separate more thorough 
treatment to be developed elsewhere. For describing manifolds in rational homotopy
theory, one transitions from using $\Q$-coefficients to  $\R$-coefficients.
However, if one considers geometry then there would be important and subtle
differences, as witnessed explicitly for instance in differential cohomology \cite{GS}. 
Consequently, this would lead to action functionals in physics taking values in $\R/\Q$,
which would require separate treatment. This together with applications, along 
the lines of \cite{top}\cite{E8}, will be discussed elsewhere. 
Note also that the description of rational 
higher connected covers in this paper should be related to the description of 
their Morava K-theory in \cite{SY}, as Morava K-theory at chromatic level zero 
is essentially rational cohomology. So some of the complementary torsion 
information not considered here is supplied in \cite{SY}.

\medskip
\noindent {\bf Notation.} We will use the notation ${\rm O}\langle n\rangle$ to denote the $(n-1)$-connected
cover of the stable orthogonal group $O$. We use $B{\rm O}\langle n+1 \rangle$ to mean 
$B({\rm O}\langle n \rangle)$, so that for our cases of interest
$BO \langle 4 \rangle =B ({\rm O}\langle 3 \rangle)=B{\rm Spin}$, $BO \langle 8 \rangle =B ({\rm O}\langle 7 \rangle)=B{\rm String}$,
$BO \langle 12 \rangle =B ({\rm O}\langle 11 \rangle)=B{\rm Fivebrane}$,
and 
$BO \langle 16 \rangle =B ({\rm O}\langle 15 \rangle)=B{\rm Ninebrane}$.

\newpage
\section{Rationalization of spaces and their approximations}


\subsection{The rational Postnikov and Whitehead towers}
\label{sectionwhite}

We start by briefly reviewing some relevant concepts in rational homotopy theory; 
see \cite{He}\cite{FHT}\cite{FOT}\cite{GM}\cite{HMR}.

\begin{definition}
\label{Def rat}
A {\rm rational space} $Y$ is a space for which   the homotopy groups 
$\pi_*(Y)$,  
or the homology groups $\tilde H_*(Y;\Z)$, or   
 $\tilde H_*(\Omega Y;\Z)$ is (equivalently) a vector space over $\Q$,
 where $\Omega Y$ is the based loop space of $Y$.
\end{definition}



\begin{definition}
\label{Def 2}
A {\rm rationalization} of a space $X$ prescribes a map $\ell_X:X\to X_\Q$ to a rational space 
$X_\Q$ such that $\ell_{X*}:\pi_*(X)\otimes\Q\to\pi_*(X_\Q)$ is an isomorphism, and
 for every map $f:X\to Y$ where $Y$ is a rational space, there is a factorization 
$$
\xymatrix@R=1.5em{
X \ar[rr]^f \ar[d]_-{\ell_X} && Y
\\
X_{\Q}
\ar @{..>}_h [urr] &&
}
$$
which is unique up to homotopy.
\end{definition}
The second property  in Definition \ref{Def 2} tells us that the space $X_\Q$ is unique up to homotopy.  Suppose $X'_\Q$ is another rational space and we have a map $\ell'_{X}:X\to X'_\Q$ satisfying the second property. Then there is a homotopy equivalence $h:X_\Q\to X'_\Q$, unique up to homotopy, such that
$$
\xymatrix@R=.5em{
&&& X_\Q \ar @/^{.5pc}/ [dd]^h 
\\
X \ar[urrr]^{\ell_X}
\ar[drrr]_{\ell_X'}
&
\\
&&&X'_\Q
\ar @/^{.5pc}/[uu]^{h^{-1}}
}
$$
is a homotopy commutative diagram. Note that an abelian group $G$ is said to be $\Q$-local 
if the map $G \otimes_\Z \Z \to G \otimes_\Z \Q$ induced by the inclusion
$\Z \hookrightarrow \Q$ is an isomorphism. 
Then Definition \ref{Def rat} can be restated as saying that
 the following holds \cite{Su2}
$$
X~ \text{is a rational space} \quad \Longleftrightarrow \quad \pi_n(X)~{\rm is}~ \Q{\rm -local} 
\quad \Longleftrightarrow \quad \widetilde{H}_n(X; \Z) ~{\rm is}~ \Q{\rm -local}\;, 
$$
for $X$ nilpotent, i.e. if its fundamental group $\pi_1$ is a nilpotent group and if $\pi_1$ acts nilpotently on the higher homotopy groups.  In particular, any simply connected space is trivially nilpotent. 
Note, however, that an extension to 
path connected spaces with general (not necessarily nilpotent) $\pi_1$ 
is possible (see \cite{FHT2}).  This allows us, for 
instance, to start our Whitehead tower (Example \ref{Ex White} below) with ${\rm SO}(n)$ or 
$B{\rm O}(n)$,  as in \cite{SS}.

\begin{example}[Eilenberg-MacLane spaces]
Consider the integral Eilenberg-MacLane space $K(\Z,n)$.  
Then the map $\iota_\Q:K(\Z,n)\to K(\Q,n)$ corresponding to the 
generator $[\iota_\Q]\in H^n(K(\Z,n),\Q)$ is a rationalization of $K(\Z,n)$.  In general, an Eilenberg-MacLane space $K(\pi, n)$ can be rationalized to 
$K(\pi \otimes \Q, n)$, as induced by the natural homomorphism 
$\pi \to \pi \otimes \Q$. 
By induction and use of the Serre spectral sequence one
can show that  
 $H^*(K(\Q, 2n); \Q)$ is a $\Q$-polynomial algebra on one generator of
degree $2n$, while $H^*(K(\Q, 2n+1); \Q)$ is a $\Q$-exterior algebra on one generator of
degree $2n+1$. Furthermore, 
$K(\Z, n) \to K(\Q, n)$ induces an isomorphism on rational 
cohomology. That is, there are isomorphisms (see \cite{GM}\cite{Mor})
\bea
H^*(K(\Z, 2n); \Q) &\cong& H^*(K(\Q, 2n); \Q) \cong \Q[\iota]\;,
\nonumber\\
H^*(K(\Z, 2n+1); \Q) &\cong& H^*(K(\Q, 2n+1); \Q) \cong \Exterior_\Q[\iota]\;,
\eea
where $\iota \in H^n(K(\Q, n); \Q)$ is the fundamental cohomology class.
 \end{example}

To construct a rationalization for a simply connected space $X$, one can use the Postnikov 
tower decomposition of $X$ and localize at each step of the tower.  
Note that one of the main properties of the rational localization functor 
is that  it preserves sequences and split suspensions, as well as commutes with 
fiber and cofiber sequences for simply connected spaces \cite{Su2}.

\begin{example}[Postnikov tower]
Recall that the Postnikov tower of $X$ (see e.g. \cite{Hat}) is a sequence of spaces and maps
	$$
	\xymatrix{
	X\ar[r] &\cdots \ar[r]& X_{(n)}\ar[r]^{\varphi_{n}} & 
	\cdots \ar[r]^{\varphi_2} & 
	X_{(1)} \ar[r]^{\varphi_1} & X_{(0)}
	}\;,
	$$
	where the map $\varphi_n$ is a fibration for every $n$,
for each space $X_{(n)}$ one has $\pi_i(X_{(n)})=0$ for $i>n$, and 
	the induced map $\varphi_{n*}:\pi_i(X_{(n)})\to\pi_i(X_{(n-1)})$ is an isomorphism for $i<n$.
%
%
%
Now assume that we have a localization $\ell_{(n-1)}:X_{(n-1)}\to X_{(n-1)\Q}$.  The principal fibration $X_{(n)} \to X_{(n-1)}$ is the pullback of the path fibration
$
PK(\pi_n(X), n+1) \to K(\pi_n(X), n+1)
$
by the Postnikov invariant $k^{n+1} \in H^{n+1}(X_{(n-1)}, \pi_n(X))$.
Now define $X_{(n)\Q}$ to be the pullback of the fibration
$
PK(\pi_n(X)\otimes \Q, n+1) \to K(\pi_n(X)\otimes \Q, n+1)
$
by the rationalized Postnikov invariant 
$$
k^{n+1}\otimes \Q \in H^{n+1}(X_{(n-1)\Q},~ \pi_n(X)\otimes \Q)\;.
$$
This then gives the localization map $\ell_{(n)}: X_{(n)} \to X_{(n)\Q}$,
completing the induction, as the base case of the induction holds because $X$ is simply connected. Defining $X_\Q$ as the inverse limit $X_\Q:=\varprojlim X_{(n)\Q}$, the rationalization of $X$ is then the induced map $\ell_X:X\to X_\Q$.
\end{example}

\medskip
We now consider the dual notion of the Whitehead tower \cite{Wh} and its
rationalization. 
 
 \begin{example}[Whitehead tower]
\label{Ex White}
  The Whitehead tower is a sequence of spaces
	$$
	\xymatrix{
	\cdots\ar[r] & 
	X\langle n\rangle\ar[r] &
	 X\langle n-1\rangle \ar[r]&
	 \cdots \ar[r] &
	  X\langle1\rangle \ar[r]&
	   X
	}\;,
	$$
where 
each space $X\langle n\rangle$ is $(n-1)$-connected, and
the map $X\langle n\rangle\to X$ induces isomorphisms on the homotopy groups $\pi_i$ for each $i\geq n$.
As with the Postnikov tower, the Whitehead tower may be constructed by induction.  Given an 
$(n-1)$-connected cover $X\langle n\rangle$ of $X$, there is a map $w_n: X\langle n\rangle\to K(\pi_{n}(X),n)$ corresponding to the generator of $H^{n}(X\langle n\rangle;\pi_{n}(X))$.  Then the space $X\langle n+1\rangle$ is constructed as the homotopy fiber of $w_{n}$. 
Similarly, one may construct the Whitehead tower over the rationalization $X_\Q$ of $X$
(see e.g. Ch. 2 in \cite{Kramer}).  By choosing $w_n\otimes \Q: X\langle n\rangle\to K(\pi_{n}(X)\otimes \Q,n)$ as the classifying map for the generator of the cohomology group 
$$
H^{n}(X\langle n\rangle_\Q;\pi_{n}(X_\Q))\cong H^{n}(X\langle n\rangle_\Q;\pi_{n}(X)\otimes \Q)\;,
$$
then the homotopy pullback $X\langle n+1\rangle_\Q$ is the rationalization of $X\langle n+1\rangle$. This also follows by induction.  Given that $\ell_{n}:X\langle n\rangle \to X \langle n \rangle_\Q$ is a localization of $X \langle n\rangle$,  we can consider the following commutative diagram
{\small 
$$
\xymatrix@=1em{
X\langle n+1\rangle\ar[dd]\ar[rr]\ar @{..>}^{\ell_{n+1}} [dr] && 
PK(\pi_n(X),n) \ar@{->}'[d][dd]\ar[dr] &
\\
& X\langle n+1\rangle_\Q\ar[rr]\ar[dd] && PK(\pi_n(X_\Q),n) \ar[dd]
\\
X\langle n\rangle\ar[dr]^{\ell_{n}}\ar@{->}'[r]^>>>>>>{w_{n}}[rr] && K(\pi_n(X),n)\ar[dr]
\\
& X\langle n\rangle_\Q\ar[rr]^-{w_{n}\otimes\Q} && ~K(\pi_n(X_\Q),n)\;.
}
$$
}
The map $\ell_{n+1}:X\langle n\rangle\to X\langle n\rangle_\Q$ exists by the universal property of pullbacks and, by the commutativity of the diagram, this map is a rational 
homotopy equivalence.
\end{example}

\begin{remark}
In the following sections we will be concerned with the rational Whitehead tower of $B{\rm O}$.  In light of this, we briefly describe one important quality of the rational Whitehead tower.  The general construction of each stage of the Whitehead tower is as the homotopy fiber of a map representing a generator of cohomology.  Consider a map $f:X\to K(A,n)$ and suppose that $f$ represents a generator of $H^n(X;A)$.  Composing $f$ with the rationalization $\iota:K(A,n)\to K(A\otimes\Q,n)$ gives a generator $f^*\iota\in H^n(A;\Q)$.  Then for any $r\in\Q$, the class $r\cdot f^*\iota$ is also a generator.  Now by the universal property of rationalization, there is a map $f_\Q:X_\Q\to K(A\otimes \Q,n)$ as well as a map $r\cdot f_\Q:X_\Q\to K(A\otimes\Q,n)$.  Denote the homotopy fibers of each of theses maps as $F(f_\Q)$ and $F(rf_\Q)$.  The key point here is to note that these spaces are homotopy equivalent.  This follows from the fact the map $r:K(\Q,n)\to K(\Q,n)$ representing $r$ times the identity induces an isomorphism on the rational cohomology of rational spaces and thus is a homotopy equivalence combined with the fact that $r\cdot f^*\iota=f^*\iota^*r$.   This is in contrast to the case where we compare $\text{hofib}(f)$ and $\text{hofib}(rf)$ for some integer $r\neq\pm1$.  These are not homotopy equivalent spaces for the reason that $r$ is not a unit in $\Z$. 
\end{remark}

\subsection{Rationalizing higher connected covers}
\label{Sec rat conn}

The groups we consider will be rationally built out of spheres.  
 Since the homotopy groups $\pi_j(S^{2m+1})$ of an odd-dimensional
  sphere $S^{2m+1}$ are finite except in top degree $j=2m+1$, 
the rationalization is $\ell_\Q: S^{2m+1} \to K(\Q, 2m+1)$
corresponding to a nontrivial class in $H^{2m+1}(S^{2m+1}; \Q)$. 
Note that we have rational equivalences 
$
S^{2m+1} \simeq_\Q K(\Z, 2m+1) \simeq_\Q K(\Q, 2m+1)
$,
i.e., the rational homotopy groups are given as 
$$
\pi_i(S^{2n+1}) \otimes \Q 
\cong 
\left\{ 
\begin{array}{ll}
\Q, & \text{if}\; i=2n+1
\\
0, & \text{otherwise}.
\end{array}
\right.
$$ 
As a  result, odd-dimensional spheres have the rational 
homotopy type of an abelian topological group,
obtained by iteratively applying the classifying space 
functor.

\medskip
The spaces we will consider are all H-spaces, and in particular H-spaces are formal.  This means that the rational homotopy type of these spaces is completely determined by their cohomology ring.  For a finite-dimensional Lie group $G$, its rational cohomology is the
same as that of a product of odd-dimensional spheres 
$$
 H^*(G; \Q) \cong H^*(\prod S^{2i-1}; \Q)\;,
$$
which implies that $G$ has the same rational homotopy type as 
that product. For our main example, 
$$
{\rm Spin}(2n) \simeq_\Q  S^{2n-1} \times S^3 \times S^7 \times \cdots \times S^{4n-5}\;, 
\qquad  \qquad
{\rm Spin}(2n+1) \simeq_\Q S^3 \times S^7 \times \cdots \times S^{4n-1}\;.
$$
 
\medskip
We can consider the question of rationalization of higher topological groups
in some generality.
%
Suppose $G$ is a topological group having the homotopy type of a CW complex. 
Rationalization commutes with products only up to homotopy. This implies that 
the rationalization of the product structure gives a map
$
\mu_\Q: G_\Q \times G_\Q \to G_\Q
$,
which may fail to be a multiplication. However, it is a group-like H-space. 
The map $\ell_\Q: G \to G_\Q$ is a homomorphism of H-spaces, i.e. the diagram 
$$
\xymatrix@R=1.5em{
G \times G \ar[rr]^\mu \ar[d]_{\ell_\Q \times \ell_\Q} && G \ar[d]^{\ell_\Q}
\\
G_\Q \times G_\Q \ar[rr]^{\mu_\Q} && G_\Q
}
$$
commutes up to homotopy, and the homomorphism is 
compatible with homotopy 
associativity. 
When $G$ is connected, 
the commutator map $\mu_\Q: G_\Q \times G_\Q \to G_\Q$  is null homotopic \cite{KSS}.

\medskip
Note that the resulting rationalizations of the connected cover groups and their classifying 
spaces in this section end up having nice models. Explicitly, they end up being products of 
rational Eilenberg-MacLane spaces. Spaces are in general twisted products of Eilenberg-MacLane 
spaces  arranged by their Postnikov decompositions. So in this case rationalization trivializes 
the Postnikov $k$-invariants, thereby untwisting the product of Eilenberg-MacLane spaces into 
a straight product. This can also be discussed on general grounds. Note that rationalization does 
not take a space outside the convenient category of CW complexes: If $X$ is a topological space 
which is a CW-complex then the rationalization $X_\Q$ is also a CW-complex.

\medskip
When $G$ is a connected topological group having the homotopy type of a finite CW complex, 
then we have that $G_\Q$ is a homotopy commutative H-space. It is homotopy equivalent as 
an H-space to a product of Eilenberg-MacLane spaces with {\it standard loop multiplication} 
\cite{KSS}.  A priori, a product of Eilenberg-MacLane spaces in general may admit many 
non-H-equivalent H-structures and thus the standard loop multiplication is not necessarily 
unique, even up to homotopy (see \cite{Cu}).  However in the rational case, as we are considering 
here, this multiplication is unique up 
to homotopy \cite{LPSS}.  

\medskip
For relatively low $k$  the connected covers $O(n) \langle k \rangle$
are defined as the based loop spaces of the corresponding classifying 
spaces in \cite{SSS3}\cite{FSSt}\cite{9brane}. It follows from the works of Kan and Milnor that every based 
loop space has the homotopy type of a topological group. In the homotopy category of connected 
CW complexes, there is an equivalence between loop spaces, topological groups,
and associative H-spaces (see \cite[Ch. 4]{Kane}). 


\begin{proposition} 
\label{prop rat grp}
In the Whitehead tower of the orthogonal group, each element in
the sequence of connected covering spaces \{${\rm String}(n)_\Q$,
${\rm Fivebrane}(n)_\Q$, ${\rm Ninebrane}(n)_\Q$, $\cdots$\}
is an abelian topological group, and this group structure is unique up to rational H-equivalence.
\end{proposition}

\theproof  (Outline)
In the context of Lie groups, one can form a product on the rational homotopy groups called the 
Samelson product.  Without defining this product, we have the following properties. From \cite{Wo}, 
if $G$ is a (possibly infinite-dimensional) connected Lie group, then the rational 
Samelson product vanishes. From \cite{LPSS}, 
for $(G,\mu)$ a connected CW homotopy-associative H-space, the following are equivalent:
	\begin{enumerate}
		\item $(G,\mu)$ has the rational H-type of an abelian topological group.
		\item $(G,\mu)$ is rationally homotopy-abelian.
		\item The Samelson bracket vanishes.
		\item There is a rational H-equivalence 
			%
		$$e:G\longrightarrow \prod_{j\geq 1} K(\pi_j(G)\otimes\Q,j)$$
			where the product of Eilenberg-MacLane spaces has the standard multiplication.
	\end{enumerate}
So, given a Lie group $G$, its Samelson product vanishes and thus $G$ is rationally homotopy-abelian. 
Moreover, for the orthogonal group as well as its connected covers, these groups are rationally equivalent 
to products of $K(\Q,n)$ spaces. Considering $\Q$ as an abelian group, it follows that $K(\Q,n)\simeq B^n\Q$ 
is also an abelian group.  Thus for groups involved in the Whitehead tower of $\Or(n)$, each space admits a rationalization by an abelian group.  Moreover, for our types of groups $G$ this rational equivalence will be multiplicative and will correspond to the unique abelian multiplication coming from the standard multiplication on Eilenberg-MacLane spaces (see \cite[Corollary 4.26 and 4.27]{LPSS}). 
\endofproof

We now consider compatibility of rationalization with taking connected covers. 
Indeed, such a problem can be studied systematically, building on classical results. 
Consider localization  $L_f$
at every prime, i.e. with respect to a map $f: \bigvee_p B\Z/p \to \ast$ with domain the 
infinite bouquet. Let $E$ be defined by the homotopy pullback diagram 
\(
\xymatrix@R=1.5em{
E \ar[rr] \ar[d]
&& X \langle n \rangle_\Q \ar[d]
\\
X \ar[rr] && X_\Q\;,
}
\label{E pull}
\)
where $X$ is a simply connected finite complex with $\pi_2(X)$ finite. 
Then  $L_fX\langle n \rangle$ has the homotopy type of $E$ \cite{Ne}. 
We now connect this to the Whitehead tower of the rationalized orthogonal group.
Since every level of that Whitehead tower is a rationalization \cite[Theorem A.1]{Dr}
this implies that at every level we will have a space which is 
at least an $H$-space. These will include all connected covers. 
Let us apply the above localization \eqref{E pull} 
to $X={\rm Spin}(n)$, which is 2-connected. 

\begin{example} [Rationalization of the 3-connected cover of ${\rm Spin}(n)$]
Let ${\rm Spin}_\tau$ be the homotopy fiber of the 
rationalization ${\rm Spin} \to {\rm Spin}_\Q$.
Starting with the fiber sequence 
$K(\Z, 2) \to {\rm String}(n)={\rm Spin}(n) \langle 7 \rangle \to {\rm Spin}(n)$,
we have $K(\Z, 2) \otimes \Q \simeq K(\Q, 2)$ as the homotopy fiber of 
$E \to {\rm Spin}(n)$ and $K(\Z, 2)_\tau$ as the homotopy fiber of 
${\rm String}(n) \to E$. 
The hypotheses imply that $E$ has the same rational homology as ${\rm Spin}(n)$,
so that $E$ is $f$-local and that $K(\Z, 2)_\tau$ is connected. The homotopy groups
of the latter are locally finite. Now the main point (see \cite{Ne}) 
 is that if $Y$ is an $f$-local space, 
then by studying the Postnikov factorization, 
every map ${\rm String}(n) \to Y$ has up to homotopy a unique extension 
to a map $E \to Y$. This identifies $L_f {\rm String}(n)$ as $E$.
\end{example}

We now consider the rationalization  of the classifying spaces of Lie groups
and their higher connected covers. 
For classifying spaces, the Borel-Hopf theorem implies that
the map $(BG)_\Q \to B(G_\Q)$ induced by the homomorphism 
$\ell_\Q: G \to G_\Q$ is a homotopy equivalence, and for connected $G$ 
we have
$$
H^*(BG; \Q)\cong \Q[y_1, \cdots, y_k]\;,
$$
where each generator $y_i$ is of even degree. 
Note that for the classifying spaces in terms of 
even generators, one cannot deduce a similar relation to spheres as above,
since the rational model for even-dimensional spheres is not free. 
It is known that the minimal model (see Sec. \ref{Sec min})
of $BG$ is evenly generated and has zero differential, i.e. $BG$ is 
rationally a product of even-dimensional Eilenberg-MacLane spaces. 
The above argument also 
applies if we have infinitely many $y_i$'s in even degrees, as in the case for
$B{\rm String}$, $B{\rm Fivebrane}$ and $B{\rm Ninebrane}$, as long as   
there are only finitely many $y_i$ in each even degree. This follows from the 
fact that two nilpotent spaces with finite Betti numbers are rationally homotopy 
equivalent if and only if they have isomorphic minimal models. 
When $G$ is connected,  the classifying space $BG$ has the rational homotopy type of 
a generalized Eilenberg-MacLane space (see \cite{FO} for another explicit description) 
and, in particular, it is rationally homotopy equivalent to a loop space \cite{KSS}.  
This implies that we can deloop as much as we desire.

\medskip
Note that the behavior of the rationalization of a Lie group $G$ is intimately 
linked to that of its classifying space $BG$.  Two compact Lie groups $G$ and $H$ 
are isomorphic  if and only if their classifying spaces $BG$ and $BH$ are homotopy 
equivalent \cite{Mo}\cite{Os}\cite{No}. 
The equivalences at the rational level are established in \cite{Mo}.
Since, rationally, taking the classifying functor only shifts the degree of the generators, 
then in this case we have that, rationally, $BG$ is a product of 
even-dimensional Eilenberg-MacLane spaces. 
Similarly in this case, taking connected covers
gives the rational models for the higher connected covers of the classifying spaces. 
Indeed, for the case of String this can be found explicitly for instance via the Serre spectral 
sequence (see \cite{SSS2} \cite{BS}). 
The higher cases work similarly. Therefore, we have 

\begin{proposition}
Given that the rational cohomology of $BG_\Q$ is $\Q[x_{2i}]$ (where $i$ can be odd or even),
the rational cohomology of the connected cover $BG \langle n \rangle_\Q$ is 
$\Q[x_{2j}]$, where $j$ runs over the subset of  values of $i$ such that $2i\geq n$. 
The same holds for the rational homotopy of the space.
\end{proposition} 

Another way to relate $G\langle n \rangle_\Q$ to $B(G\langle n \rangle)_\Q=BG \langle n+1 \rangle_\Q$ is  as follows. Rationalization is a weak rational equivalence $X \to L_\Q(X)$, where $L_\Q$ is the $\Q$-localization functor. 
Suspension and looping preserve the rationality of the 
spaces involved. More precisely, for any simply connected space $X$, the loop space
of its rationalization is 
$\Omega X_\Q \simeq \prod_\alpha K(\Q, m_\alpha)$ for various values of $m_\alpha$ \cite{Str}.
Note that in one approach
to the String group, it is taken to be the loop space of its classifying space, with the latter
constructed first (see \cite{SSS3}).  Indeed, for $Y=BG\langle n \rangle$ 
a pointed topological space, the assignment $X \mapsto [X, \Omega B G \langle n \rangle]_*
:={\rm Top}_*(X, \Omega B G \langle n \rangle)/\simeq_*$ defines a functor 
$[-, \Omega B G \langle n \rangle]_*={\rm Top}_*^{\rm op} \to {\rm Grp}$ 
from the category of pointed topological spaces to the category of groups, where
the multiplication in $[X, \Omega B G \langle n \rangle]_*$ is pointwise 
(see \cite[Thm. 1.2.5]{Pic} for the corresponding classical statement).
Starting with $X= B{\rm String}$ and rationalizing $X_\Q\simeq (B{\rm String})_\Q$, then taking
 the loop space leads to
$$
{\rm String}_\Q \simeq \Omega (B{\rm String})_\Q \simeq \prod_\alpha K(\Q, m_\alpha)\;.
$$
The values of $m_\alpha$ can be deduced from another approach, namely starting with 
the Spin group, taking connected covers and then rationalizing, as we did earlier. 
A similar treatment for Fivebrane and Ninebrane can be established in parallel. 
Overall, we can arrive at the result 
that the rational models for ${\rm Spin} (n) \langle k \rangle_\Q$ 
and $B{\rm Spin} (n)\langle k \rangle_\Q$ are indeed given as products 
of rational Eilenberg-MacLane spaces. 

\medskip
We next consider another 
way of obtaining the rational homotopy type from the rational 
cohomology ring and vice versa.

\subsection{Minimal Models}
\label{Sec min}

Now we present minimal models (see \cite{FHT}\cite{FHT2}\cite{FOT}\cite{GM}\cite{BG}) 
for our connected covers in a straightforward manner.
Given a space $X$, one can assign the complex of piecewise linear differential forms, $A^*_{PL}(X)$ (see \cite{GM} for a definition).  This assignment defines a functor from the category of topological spaces to the category of CDGAs over the rationals.  A minimal model for a CDGA $(A,d)$ is given by a quasi-isomorphism $\varphi:(\Exterior V,d_V)\to(A,d_A)$ where $\Exterior V$ is freely generated and the image of $d_V$ is contained in the set of decomposable elements of $\Exterior V$. For a space $X$, a minimal model for $X$ is defined to be a minimal model for $A^*_{PL}(X)$.  

\medskip
Let $E\xrightarrow{p}B$ be a quasi-nilpotent fibration
with fiber $F$, let $f:B'\to B$ be a map of base spaces, 
and let $E'\xrightarrow{p'}B'$ denote the pullback of $p$ along $f$.  Letting $(\Exterior V,d_V)\to A_{PL}(B)$ be a minimal model for $B$ then the relative minimal model for the fibration $p$ is given by $(\Exterior V,d_V)\to(\Exterior V\otimes \Exterior W, d)\to(\Exterior W,\bar d)$ where $(\Exterior V,d_V)\to(\Exterior V\otimes\Exterior W,d)$ is the relative minimal model for $p$ as a map and $(\Exterior W,\bar d)$ is formed as the quotient $(\Exterior V\otimes\Exterior W,d)/(\Exterior^+(V)\otimes\Exterior W)$. Let $\phi:(\Exterior V,d)\to(\Exterior V',d')$ denote the relative minimal model for $f$.  The following standard result from rational homotopy theory provides a recipe for constructing a minimal model of a pullback (see \cite{FHT} \cite{FOT}).
The relative minimal model 
	$$
	\xymatrix{
	(\Exterior V',d')\ar[r] &
	(\Exterior V'\otimes \Exterior W,D)
	\ar[r]^-{\rho} &
	(\Exterior W,\bar D)
	}
	$$
	is the relative minimal model for the pullback fibration $p'$. The CDGA $(\Exterior V'\otimes\Exterior W,D)$ is defined as $(\Exterior V',d')\otimes_{\Exterior V}(\Exterior V\otimes\Exterior W,D)$.  Here $D$ is defined by $D(w)=(\phi\otimes1)(d'w)$ where $\phi\otimes1:\Exterior V\otimes\Exterior W\to\Exterior V'\otimes\Exterior W$, and 
$(\Exterior W,\overline{D})$ is obtained by a quotient.


\medskip
Consider the rational Whitehead tower of $B\Or(n)$.  
Recall  (Example \ref{Ex White}) that this tower can be constructed as a system of pullbacks where at each step of the tower, the space $B(\Or(n)\langle 4k+3\rangle)$ is formed via the pullback of the pathspace fibration $PK(\Q,4k)\to K(\Q,4k)$ along a map
 $p_{k}:B(\Or(n)\langle 4k-1\rangle)\to K(\Q,4k)$ which can be thought of as a rationalization of the  $k$-th 
Pontrjagin class.   

\medskip

The relative minimal model for the pathspace fibration of a space $X$ is given by 
$(\Exterior V,d)\to(\Exterior V\otimes\Exterior sV,d)\to(\Exterior sV,\bar d)$ 
where an element $sv\in\Exterior sV$ has degree $|sv|=|v|-1$.  As the minimal 
model for the Eilenberg-MacLane space $K(\Q,4k)$ is given by $(\Exterior(y_{4k}),0)$, 
then it follows that the minimal model for the pathspace fibration is given by 
$(\Exterior(y_{4k}),0)\to(\Exterior(y_{4k},sy_{4k}),d)\to(\Exterior(sy_{4k}),0)$ 
where the differential $d$ is given by $d(y_{4k})=0$ and $d(sy_{4k})=y_{4k}$.  
Then we immediately have the following.

\begin{proposition} 
A relative minimal model for the fibration $B({\rm O}(n)\langle 4p+3\rangle)
 \to B({\rm O}(n)\langle 4p-1\rangle)$, i.e. of the fibration 
 $B {\rm O}(n)\langle 4p+4 \rangle \to B {\rm O}(n) \langle 4p \rangle$, is given by 
	\begin{enumerate}
	\item[(a)] For $n=2k$:
	$$
	\xymatrix{
	\big(\Exterior (x_{4p},x_{4p+4},\ldots,x_{4k-4},\chi_{2k}),~0\big)
	\ar[r] & \big(\Exterior (x_{4p},\ldots,x_{4k-4},\chi_{2k},sy_{4p}),~d\big)
	\ar[r] & \big(\Exterior(sy_{4p}),~0\big)\,,
	}
	$$
where $d(sy_{4p})=x_{4p}$ and 0 otherwise.\\
	\item[(b)] For $n=2k+1$:
	$$
	\xymatrix{
	\big(\Exterior (x_{4p},x_{4p+4},\ldots,x_{4k}),~0\big)
	\ar[r] & \big(\Exterior (x_{4p},\ldots,x_{4k},sy_{4p}),~d\big)
	\ar[r] & \big(\Exterior(sy_{4p}),~0\big)\,,
	}
	$$
where $d(sy_{4p})=x_{4p}$ and 0 otherwise.
\end{enumerate}
\end{proposition}
Note that the element $\chi_{2k}$  corresponds to the 
Euler class for even $n$. 
We now consider explicitly the cases $B{\rm Fivebrane}(n):=B{\rm O}(n)\langle 9 \rangle$ 
and $B{\rm Ninebrane}(n):=B{\rm O}(n)\langle 13 \rangle$.

\begin{example}
{\bf (i)} In the case of  the fibration $B\Five(n)\to B\String(n)$, it follows that for $n$ odd the minimal model is 
	$$
	\xymatrix{
	(\Exterior(x_8,x_{12},\ldots,x_{4k}),~0)\ar[r] & 
	(\Exterior(x_8,\cdots,x_{4k},sy_{8}),~d)\ar[r] &
	(\Exterior(sy_8),~0)\,,
	}
	$$
where $x_8$ is the element corresponding to the second Pontrjagin class.  
\item {\bf (ii)} In the case of  the fibration $B\Nine(n)\to B\Five(n)$ and $n$ odd, 
the minimal model is 
	$$
	\xymatrix{
	(\Exterior(x_{12},\ldots,x_{4n+4}),~0)
	\ar[r] &
	(\Exterior(x_{12},\cdots,x_{4n+4}, sy_{12}),~d)\ar[r]&
	(\Exterior(sy_{12}),~0)\,,
	}
	$$
where $x_{12}$ is the element corresponding to the third Pontrjagin class.  These can be modified appropriately by adding the Euler class in the case when $n$ is even.

\end{example} 

Note that, as was mentioned in Sec. \ref{Sec rat conn}, all the spaces we consider are H-spaces and thus formal. Nevertheless, we believe that
 it is interesting to explicitly present the minimal models as above.

\subsection{Rank vs. connectivity degree}
\label{Sec rank}

In this section we will highlight how the rank $n$ of the Spin group
${\rm Spin}(n)$ will have an effect on the corresponding $k$-connected
cover. We start with identifying the minimal rank so that the resulting 
rationalizations are not trivial, after which we consider the indefinite case 
${\rm Spin}(p, q)$.

\begin{example}[The unstable case: ${\rm String}(3)$]
\label{String3}
The case $n=3$ is special. From the point of view of 
classifying spaces, the generator $Q_1=\tfrac{1}{2}p_1$ in 
$H^4(B{\rm Spin}(3); \Z)$ is further divisible by 2, or 
the first Pontrjagin class pulled back from $B{\rm SO}(3)$ 
via the covering map ${\rm Spin}(3) \to {\rm O}(3)$ is divisible by 4.
This has interesting consequences that we will not pursue here
(see \cite{Red}\cite{tcu}).
From the rational point of view, however, the story is different. Consider
the identification ${\rm Spin}(3) \cong S^3$. Now forming 
${\rm String}(3)$ is equivalent to forming the 3-connected cover
$S^3\langle 4 \rangle$ of the 3-sphere. This latter space is known 
to be torsion. This is essentially due to Serre's result that 
$\pi_j(S^3)$ is finite for $j>3$. 
 Indeed forming the fibration 
$
{\rm String}(3) \to {\rm Spin}(3) \to K(\Z, 3)
$
and rationalizing, we consider the fibration 
$
{\rm Spin}_\tau(3) \to {\rm Spin}(3) \to {\rm Spin}(3)_\Q
$,
where the leftmost term is the homotopy fiber to be determined. 
This is homotopy equivalent to the fibration
$
S^3_\tau \to S^3 \to S^3_\Q
$,
where $S^3_\Q \simeq K(\Q, 3) \simeq M(\Q, 3)$ 
and  the homotopy fiber $S^3_\tau \simeq M(\Q/\Z, 2)$
is the Moore space for the quotient $\Q/\Z$. 
\end{example}
 
 
In order to generalize to higher connected covers, we note that in the case of $n=3$, 
we have $\dim({\rm Spin}(n))=3$, and the generator of $\pi_3(\Spin(3))$ corresponds to the
 fundamental class in $H^3(\Spin(3);\Z)$.  In the general case the dimension of Spin($n$) is 
 $d=\frac{1}{2}n(n-1)$.  Now, while it is true that 
 $\Spin(n)\langle\frac{1}{2}n(n-1)\rangle$ are torsion spaces and thus rationally trivial, 
 we can state a sharper result.

\begin{proposition} 
The $(k-1)$-connected cover of rank $n\geq 2$, ${\rm Spin}(n) \langle k \rangle_\Q$,
 is homotopy trivial for $k\geq4\cdot\lfloor\tfrac{n-1}{2}\rfloor$. 
\end{proposition}
\proof
It is a classic result that the rational cohomology $H^*(\Spin(n);\Q)$ is isomorphic to the exterior algebra 
$\Exterior_\Q(x_3,x_7,\ldots,x_{4i-1})$ when $n$ is odd and to $\Exterior_\Q(x_3,x_7,\ldots,x_{4i-1},y_{n-1})$ 
for $n$ even, where $i=\lfloor\tfrac{n-1}{2}\rfloor$.  Moreover it is clear that these algebras, equipped 
with the zero differential, describe a minimal model for $\Spin(n)$
(see Sec. \ref{Sec min}).  From rational homotopy we know that $V^n\cong\text{Hom}(\pi_n(\Spin(n)),\Q)$ 
where $V^n$ denotes the vector space of degree $n$ generators.  Thus it follows that the non-torsion generator 
of highest degree is in degree $4\cdot\lfloor\tfrac{n-1}{2}\rfloor-1$ and, upon killing this homotopy class
 in the Whitehead tower, the resulting space is pure torsion.

\vspace{-5mm}
\endofproof

Thus, for example, we have that ${\rm Fivebrane}(n):={\rm Spin}(n)\langle 8 \rangle$ is 
a torsion space for $n\leq6$, and ${\rm Ninebrane}(n):={\rm Spin}(n)\langle 12 \rangle$ 
is a torsion space for $n\leq 8$.

\paragraph{\bf The indefinite case ${\rm Spin}(p, q)\langle k \rangle$}
Note that we can consider rationalization of 
higher structures in the indefinite signature, i.e. by taking connected
covers of the semi-orthogonal $SO(p, q)$, prominent in semi-Riemannian geometry. 
For degree 3, i.e. for ${\rm String}(p, q)$, these are 
characterized in \cite{SS}. The homotopy groups encountered there are complicated, 
but upon rationalizing the problem becomes much more tractable: If 
$p, q <3$ then the problem is trivial. If we have $p=q=3$ then 
we have two copies of the trivialization of the 3-connected cover
of $S^3$, which is pure torsion (as in Example \ref{String3}). 
For $p=q=4$, we have four copies of $S^3\langle 4 \rangle$. 
So far the rationalization of all these cases is trivial. 
Once we reach  $p, q \geq 5$ then we have two copies of the nontrivial  problem, i.e.,
 a rationalization of ${\rm String}(p)$ and of ${\rm String}(q)$.  The cases when 
 $p \neq q$ can be dealt with similarly.

\medskip
Given the discussions in previous sections, the descriptions of the higher 
connected covers ${\rm Fivebrane}(p, q)_\Q$ and ${\rm Ninebrane}(p, q)$ 
rationally will follow
analogously. Note that this is in stark contrast with the calculations 
in \cite{SS} where the various torsion groups arising notoriously in the 
indefinite case made the extension to Fivebrane and Ninebrane
not immediately possible. 

\begin{proposition}
Let ${\rm Spin} (p, q)\langle k \rangle$ denote the $(k-1)$-connected cover of 
the indefinite Spin group ${\rm Spin}(p, q)$. 
\item {\bf (i)} The rationalization takes
the form 
$$
{\rm Spin} (p, q)\langle k \rangle_\Q \simeq 
{\rm Spin} (p)\langle k \rangle_\Q
\times {\rm Spin} (q)\langle k \rangle_\Q 
\qquad 
{\rm if}~ k<    4\cdot\lfloor\tfrac{{\rm min}\{p, q\}-1}{2}\rfloor\;.
$$ 
\item {\bf (ii)} This rationalization is trivial, i.e. 
the resulting spaces are torsion spaces, when 
\begin{itemize}
\item $p=q$ such that 
$k\geq4\cdot\lfloor\tfrac{p-1}{2}\rfloor=4\cdot\lfloor\tfrac{q-1}{2}\rfloor$, or
\item $p\neq q$ such that $k \geq 4\lfloor p-1 \rfloor$, $k \geq \lfloor q-1 \rfloor$.  
\end{itemize}
\end{proposition}

Intermediate cases can arise. 
For example, for $\op{Fivebrane}(7, q)$ the first factor 
$\op{Fivebrane}(7)$ is a torsion space, while the second factor 
$\op{Fivebrane}(q)$ has a nontrivial rationalization for $q>7$.

\section{Higher tangential structures}

\subsection{Rational Structures}
 
Let $G$ be a simply connected topological group.  
For a principal $G$-bundle $P\to M$ and a homomorphism $\rho:H\to G$ of topological groups, one says that the structure group of $P$ lifts from $G$ to $H$ if there is a principal $H$-bundle $Q\to M$ and a bundle isomorphism $Q\times_\rho G\cong P$ over $M$, and any principal $H$-bundle satisfying this property is an $H$-structure for $P$. Two $H$-structures are isomorphic if there is a bundle isomorphism between them.  From a homotopy theoretic perspective, we can associate to any $G$-bundle $P\to M$ a classifying map $f:M\to BG$.  The  homomorphism $\rho:H\to G$ induces a map of classifying spaces, $B\rho:BH\to BG$ and the associated $G$-bundle, $EH\times_\rho G\to BH$, is classified by the map $B\rho$.  Thus a lifting of the classifying map along $B\rho$
	$$
\xymatrix@R=1.5em{
& &BH \ar[d]^{B\rho}
\\
M \ar@{..>}[urr]^{\tilde f}
\ar[rr]^{f}
& &BG.
}
$$
corresponds to a lifting of the structure group from $G$ to $H$ as $(\tilde f^*EH)\times_\rho G\cong \tilde f^*B\rho^*EG\cong f^*EG$.  In fact given a principal $G$-bundle $P\to M$ and fixing a choice of classifying map, then there is a one-to-one correspondence between homotopy classes of lifts of the classifying map along $B\rho$ and isomorphism classes of principal $H$-bundles which represent lifts of the structure group from $G$ to $H$.  We will take the homotopy theoretic perspective in understanding lifts of the structure group and use the following definition.
\begin{definition}
Let $P\to M$ be a principal $G$-bundle. 
\begin{enumerate}
	\item An $H$ structure on $P$ is a lift of the classifying map from $BG$ to $BH$.
	\item Two $H$ structures $\tilde f,\tilde f'$ are isomorphic if there exists a homotopy $\mathcal{H}:[0,1]\times M\to BH$ such that
		$$\mathcal{H}(0,x)=\tilde f(x),\quad\mathcal{H}(1,x)=\tilde f'(x),\quad\text{and}\quad B\rho\circ \mathcal{H}(t,x)=f(x),\;\;\forall t\in[0,1].$$
\end{enumerate}
\end{definition}
\medskip
In \cite{Red}, Redden studies the general case for liftings of the structure group where $BH$ is the homotopy fiber of a map $\lambda:BG\to K(A,k)$ for an abelian group $A$ and such that the group $G$ is $(k-2)$-connected.  For the purposes of this paper, we focus on topological groups arising in the Whitehead tower for ${\rm O}(n)$ and thus the cases where $G={\rm O}(n)\langle k-1\rangle$ and $H={\rm O}(n)\langle k\rangle$. We further remark that unless otherwise stated, we assume that ${\rm O}(n)\langle k-1\rangle$ is in the stable range for $n$ and will thus drop the index $(n)$ from here on out.  Exploiting the connectivity of ${\rm O}\langle k-1\rangle$ we have $H^{k}(B{\rm O}\langle k\rangle;\pi_{k-1}({\rm O}))\cong\pi_{k-1}({\rm O})$.  For example, when $k=4$ we have $B{\rm O}\langle 4\rangle=B{\rm Spin}$ and $H^4(B\Spin;\pi_3({\rm O}))\cong\ZZ$. Combining this isomorphism with Brown's representability theorem, the generator $\theta_{k}\in\pi_{k-1}({\rm O})$ corresponds to a map $\theta_{k}:B{\rm O}\langle k\rangle \to K(\pi_{k-1}({\rm O}),k)$. For a principal ${\rm O}\langle k-1\rangle$-bundle, we denote $\theta_k(P):=f^*\theta_k$ and this class represents the obstruction for $P$ to admit an ${\rm O}\langle k\rangle$-structure.  

\medskip
Using the loop space functor, there is a morphism $[X,K(A,k)]\to[\Omega X,\Omega K(A,k)]$ which corresponds to a morphism $H^*(X;A)\to H^{*-1}(\Omega X;A)$.  Setting $X=BG$ and identifying $\Omega BG\simeq G$, we obtain the transgression map $\tau:H^*(BG;A)\to H^{*-1}(G;A)$ for any group $G$.  In fact, if the group $G$ is $(k-2)$-connected, then $\tau$ is the right inverse for $d_k$, i.e. $d_k(\tau)=\text{Id}$, where $d_k$ is the cohomological transgression arising from the $k$th page of the Serre spectral sequence of a fibration.  
Specializing to our case, we combine several results from Section 2 of \cite{Red} into the following.

\begin{proposition} Let $P\to M$ be an ${\rm O}\langle k-1\rangle$-bundle. Then 
\begin{enumerate}\label{torsorprop}
\item $P$ admits an ${\rm O}\langle k\rangle$-structure if  and only if $f^*\theta_k=0$. 
\item There is a one-to-one correspondence between homotopy classes of ${\rm O}\langle k\rangle$-structures and cohomology classes $\gamma\in H^{k-1}(P;\pi_{k-1}({\rm O}))$ such that $\iota_x^*\gamma=\tau\theta_k(P)$.
\item The set of ${\rm O}\langle k\rangle$-structures up to homotopy is an $H^{k-1}(M;\pi_{k-1}({\rm O}))$-torsor.
\end{enumerate}
\end{proposition}
We provide a sketch of the proof here.  A full and much more detailed proof of this proposition can be found in \cite{Red}.  We note that the first statement follows from the fact that $B{\rm O}\langle k+1\rangle$ can be realized as the homotopy fiber of $\theta_k:B{\rm O}\langle k\rangle\to K(\pi_{k-1}({\rm O}),k)$ and thus a lift of the classifying map $f$ exists if and only if $\theta_k\circ f\simeq*$. The second statement requires somewhat more detail.  However we note that the map from the homotopy classes of ${\rm O}\langle k\rangle$-structures to cohomology on the total space $P$ is given by using the contractibility of $E{\rm O}\langle k\rangle$.  A lift of the classifying map $f$ to $B{\rm O}\langle k+1\rangle$ induces a map on $P$ to a fiber of $B\rho^*E{\rm O}\langle k\rangle$, which by construction has the homotopy type of $K(\pi_{k-1}({\rm O}),k-1)$, and thus defines a cohomology class in $H^{k-1}(P;\pi_{k-1}({\rm O}))$. For the third statement, we can consider the Serre spectral sequence corresponding to the principal bundle $P\to M$.  Utilizing the connectivity of the fiber, there is an exact sequence 
	$$
	\xymatrix{
	0\ar[r] &
	 H^{k-1}(M;\pi_{k-1}({\rm O}))\ar[r]^-{\pi^*}&
	 H^{k-1}(P;\pi_{k-1}({\rm O})) \ar[r]^-{\iota_x^*}&
	 H^{k-1}(P_x;\pi_{k-1}({\rm O}))\ar[r]^-{d_k}&
	 H^k(M;\pi_{k-1}({\rm O}))
	}\;,
	$$
where $\iota_x^*$ denotes the pullback along the inclusion of the fiber over $x$ and $d_k$ is the differential arising from the $k$th page of the Serre spectral sequence. 
\medskip

The result of this proposition is that one can classify ${\rm O}\langle k\rangle$-structures by certain cohomology classes in the total space.  This leads us to the following definition.

\begin{definition}
For an ${\rm O}\langle k-1\rangle$-bundle $P\to M$, an {\rm ${\rm O}\langle k\rangle$ class} is a cohomology class $\gamma\in H^{k-1}(P;\pi_{k-1}(\rm O))$ such that $\iota_x^*\gamma=
\tau\theta_k(P)$ for each fiber inclusion $\iota_x:{\rm O}\langle k-1\rangle\to P$.
\end{definition} 

Considering the process of rationalization and our discussion surrounding the rational Whitehead tower in Sec. \ref{sectionwhite}, we can construct a nice parallel to this story of 
${\rm O}\langle k\rangle$-structures by considering the rational Whitehead tower over $B{\rm O}$.  Given a homomorphism of groups $\rho:H\to G$, let $\rho_\Q:H_\Q\to G_\Q$ represent the rationalization of $\rho$.  This in turn induces a morphism of classifying spaces $B\rho_\Q:BH_\Q\to BG_\Q$. 
Note that, as we have seen in Sec. \ref{Sec rat conn},
we can think about $BG_\Q$ equally via either delooping of H-spaces or via classifying spaces of groups.   
On the other hand, if we have a principal $G$-bundle, then the composition of the classifying map $f:M\to BG$ with the rationalization of $BG$, $\ell_{BG}:BG\to BG_\Q$, gives a map $f_\Q:M\to BG_\Q$.  Note that as $B(G_\Q)\simeq (BG)_\Q$ then, in the context of classifying spaces, there is no ambiguity in using the notation $BG_\Q$. In order to pursue the analogy with $G$-structures as above, we begin with the following definition.

\begin{definition}
Given a principal $G$-bundle $P\to M$, a {\rm rational $H$-structure} on $P$ is given by a lift of the classifying map $f_\Q:M\to BG_\Q$ along the map $B\rho_\Q:BH_\Q\to BG_\Q$.
\end{definition}

For our purposes, we specialize this to the case of ${\rm O}\langle k \rangle$-structures. 

\begin{example}[Rational Whitehead tower of ${\rm O}$]
\label{Ex ROn}
 Consider the rational Whitehead tower corresponding to the classifying space $B{\rm O}$.  As the rationalization induces an isomorphism on homotopy groups tensored with $\Q$ and, since the only non-torsion homotopy groups $\pi_i({\rm O})$ occur when $i=4k-1$,  the 
rational Whitehead tower looks as follows

\vspace{-7mm}
{\small
$$
\xymatrix@=1.6em{
 &\vdots \ar[d]& &\\
K(\Q,11)\ar[r] &B{\rm O}\langle 16 \rangle_\Q=B{\rm Ninebrane}_\Q \ar[d] &&\\
K(\Q,7)\ar[r] &B{\rm O}\langle 12 \rangle_\Q =B{\rm Fivebrane}_\Q 
\ar[d]\ar[rr]^-{p_3^\Q}&&K(\Q,12)\\
K(\Q,3)\ar[r]&B{\rm O}\langle 8 \rangle_\Q=B{\rm String}_\Q 
\ar[rr]^-{p_2^\Q}\ar[d]&&K(\Q,8)\\
 &B{\rm O}\langle 4 \rangle_\Q=B{\rm Spin}_\Q 
 \ar[rr]^-{p_1^\Q}&&K(\Q,4)\;.
}
$$
}

\noindent The first two homotopy groups of ${\rm O}$ are torsion, so that 
$B{\rm O}_\Q \simeq B{\rm SO}_\Q \cong B{\rm Spin}_\Q$. Similarly, in the 
next period of the real Bott periodicity,  
the two homotopy groups of ${\rm O}$ in degrees 9 and 10 are torsion, so that
(in the notation of \cite{9brane}, see also Sec. \ref{Sec var9})
$B2{\rm Orient}_\Q \simeq B2{\rm Spin}_\Q \simeq B{\rm Ninebrane}_\Q$. Note that
the obstructions are given a priori by fractions of the indicated Pontrjagin classes. However,
since we are working rationally, these are equivalent to the bare classes. 
\end{example} 

Once we rationalize, our structures connect to classical constructions.

\begin{example} 
{\bf (i)}  (Rational String structures). Let $G$ be a compact Lie group and $H$ a closed subgroup. 
Then the coset space $G/H$ has vanishing $p_1^\Q$ provided 
the Killing form of $\frak{g}={\rm Lie}(G)$ restricts to a multiple of the 
Killing form of $\frak{h}={\rm Lie}(H)$ \cite{Back}. Consequently, 
the same holds if the subgroup $H$ is simple \cite{Back}\cite{Sing}.

 \item {\bf (ii)} (Rational Fivebrane structures). Manifolds admitting differentiable action of the groups 
 $G={\rm SU}(n)$, ${\rm SO}(n)$, or 
${\rm Sp}(n)$ with vanishing $p_1^\Q$ and $p_2^\Q$ 
are studied in \cite{HH1}\cite{Gr}.

\item {\bf (iii)} (Rational Ninebrane structures). In fact, in \cite[Prop. 3.2]{Sing} a criterion for 
constructing quotients $G/H$ with $p_1^\Q=p_2^\Q=p_3^\Q$ 
is given in terms of the positive roots $\{\beta_i\}_{i=1}^s$ of $H$. 
This is the requirement that the sum $\sum_{i=1}^s \beta_i^{2m}$ is contained in the ideal
$I$ of the ring $H^*(BH; \Q)$ generated by
$q^*(\widetilde{H}^*(BG; \Q))$ for $1 \leq m \leq 3$.

\item {\bf (iv)} (Rational Ninebrane structures). 
Requiring $p_1^\Q=p_2^\Q=p_3^\Q=0$ for manifolds of dimension at most twelve 
is equivalent to these being rationally parallelizable. This happens for $G/H$ with
$H$ locally isomorphic to $SU(2)$ (see \cite[Cor 3.3]{Sing}). 

\end{example}

\subsection{Variations on rational Fivebrane classes}
\label{Sec var5}

We start by discussing in detail the $k=8$ case, i.e. when we have a 
principal $\String$-bundle and wish to investigate when it admits a rational $\Five$ structure.  As the space $\String$ is 6-connected with $\pi_7(\String)\cong\pi_7(O)\cong\Z$, it follows from the Hurewicz and Universal Coefficients Theorem that $H^7(\String;\Q)\cong\Q$. Tracing what it means for a manifold $M$ to have a Fivebrane structure, we make the following definitions.

\begin{definition}
A {\rm rational Fivebrane structure} is a lift of the 
 $\String$-principal bundle $\pi_{{}_\String}:P\to M$
 to the homotopy fiber $\text{hofib}(\tfrac{1}{6}p_2^\Q)$ 
 of the rationalized classifying map $f:M\to B\String_\Q$. 
\end{definition}

Note here that we have chosen to study the homotopy fiber of a representative of $\tfrac{1}{6}p_{2}^\Q$.  As above (see Example \ref{Ex ROn}), we could have chosen to study the homotopy fiber of $p_{2}^\Q$ or even $rp_{2}^\Q$ for any $r\in\Q$, as the resulting classifying spaces are all homotopy equivalent.  The only discrepancy here will be that the rational Fivebrane structures will differ by an homotopy equivalence.  Thus, up to isomorphism, these structures are the same.    

%
%
\medskip
We now refer  to the discussion just before Prop. \ref{torsorprop}.
Setting $a_7\in H^7(\String;\Q)$ to be the generator given by $a_7:=\tau(\psix^\Q)$, we make the following definition.


\begin{definition}
A {\rm rational Fivebrane structure class} is a cohomology class $\mathcal{F}\in H^7(P;\Q)$ such that $\iota_x^*\mathcal{F}=a_7\in H^7(\String;\Q)$ for each fiber inclusion $\iota_x:\String\to P$.
\label{Def 5}
\end{definition}

As with the integral case \cite{SSS2}, it follows that these rational Fivebrane structure classes form a torsor for $H^7(M;\Q).$  Furthermore, the case of finite rank can be treated similarly, taking into account the discussions in Sec. \ref{Sec rank}.

\medskip
At this stage, our goal is to describe higher structures (beyond Spin) using Spin structures to the extent of which it is possible.  We will do this here for Fivebrane and in the next section for Ninebrane structures.  Given a principal String-bundle $\pi_{\String}:P\to M$, there is an underlying principal Spin-bundle $\pi_{\Spin}:Q\to M$ which fits into the following commutative diagram
	$$\xymatrix@R=1em@C=4em{\String\ar[r]\ar[dd]_{\mu_0}&P\ar[rd]^{\pi_{\String}}\ar[dd]_{\mu}&\\
				&&M\\
				\Spin\ar[r]&Q\ar[ur]_{\pi_{\Spin}}&}
	$$
where the homomorphism $\mu_0:\String\to\Spin$ has fiber a $K(\ZZ,2)$ and the bundle map $\mu$ is $\mu_0$ equivariant.  For the homomorphism $\mu_0$, we have the following useful fact when considering rational cohomology.
\begin{lemma}
\label{spin5iso} The map $\mu_0:\String\to \Spin$ induces an isomorphism $\mu_0^*:H^7(\Spin;\Q)\xrightarrow{\cong}H^7(\String;\Q)$.
\end{lemma}

We have defined rational Fivebrane classes solely as any class in $H^7(P;\Q)$ which restricts to a certain generator in $H^7(\String;\Q)$.   We make two notes regarding this.  The first is that the transgression map is invariant under rationalization.  Thus if we have a generator 
$\tfrac{1}{6}p_2\in H^8(B\String;\Z)$, then $(\frac{1}{6}p_2)^\Q$ is a generator for $H^8(B\String;\Q)$.  More importantly, we have $\tau((\tfrac{1}{6}p_2)^\Q)=(\tau(\psix))^\Q$.  The consequence of this is the following.

\begin{lemma}
The rationalization of any Fivebrane class is a rational Fivebrane class.
\end{lemma}
\proof
Every Fivebrane class $\mathcal{F}\in H^7(P;\Z)$ satisfies $\iota_x^*\mathcal{F}=\tau(\psix)$.  Then,
 by naturality of rationalization and what we noted above, the rational class 
$\mathcal{F}_\Q$ satisfies
	$$
	\iota_x^*\mathcal{F}_\Q=(\iota_x^*\mathcal{F})_\Q=(\tau(\psix))_\Q=a_7\;.
	$$
Hence $\mathcal{F}_\Q$ is a rational Fivebrane class.
\endofproof

Thus for any ordinary Fivebrane class, there is a corresponding rational Fivebrane class.  The second thing we note is that with the isomorphism from Lemma \ref{spin5iso}, we can define a generator of $H^7(\Spin;\Q)$ as $(\rho^*)^{-1}(\tau(\psix))$.  For simplicity we will denote this class as $\tilde a_7$.  We will also set $a_3:=\tau((\tfrac{1}{2}p_1)^\Q)\in H^3(\Spin;\Q)$.  Consequently, by considering the underlying Spin bundle for our String bundle, we can define classes here similar to how Fivebrane classes are defined cohomologically.

\medskip
To that end, let $\pi_\Spin:Q\to M$ denote the underlying $\Spin$ bundle.
\begin{definition}
A {\rm rational Spin-Fivebrane class} is a cohomology class $\mathcal{F}_\Q$ in  $H^7(Q;\Q)$ such that $\iota_x^*\mathcal{F}_\Q=\tilde a_7\in H^7(\Spin;\Q)$ for each $x\in M$.
\label{Def Spin5}
\end{definition}

The main question we pursue now is how the two definitions,
Def. \ref{Def 5} and Def. \ref{Def Spin5}, are related.  It is not too difficult to show that every 
rational Spin-Fivebrane class gets mapped by $\mu^*$ to a rational Fivebrane class; however we can say more, still for String bundles. 

\begin{theorem} Let $\pi_{\String}:P\to M$ be a principal $\String$-bundle and let $\pi_{\Spin}:Q\to M$ be its underlying principal Spin-bundle.
\label{rfiv} {\bf (i)} For every rational Spin-Fivebrane class 
$\mathcal{F}\in H^7(Q;\Q)$, the pullback $\rho^*\mathcal{F}$ is a rational Fivebrane class.
		\item {\bf (ii)} For any rational $\Five$ class $\mathcal{F}\in H^7(P;\Q)$ there is
		 a Spin-Fivebrane class  $\tilde{\mathcal{F}} \in H^7(Q;\Q)$ such that 
		$\mu^*\tilde{\mathcal{F}}=\mathcal{F}$.
		\item {\bf (iii)} Two classes $\mathcal{F}$, $\mathcal{F}'\in H^7(Q;\Q)$ will give the same rational Fivebrane class if $\mathcal{F}-\mathcal{F}'=\mathcal{S}\cdot\pi_\Spin^*\phi_4$ where $\mathcal{S}\in H^3(Q;\Q)$ is the String structure class and 
		$\phi_4\in H^4(M;\Q)$ is a rational cohomology class.
\end{theorem}
\proof
The main ingredient that will be used in the proof is the corresponding Serre spectral sequences for the bundles $Q$ and $P$ along with the spectral sequences for the universal bundles over the classifying spaces $B\Spin$ and $B\String$.  Let $f_{\Spin}$ and $f_\String$ be the classifying map of $Q$ and $P$ respectively.
The second page of the spectral sequence for $Q$ is as follows
\begin{center}
\includegraphics[width=200px]{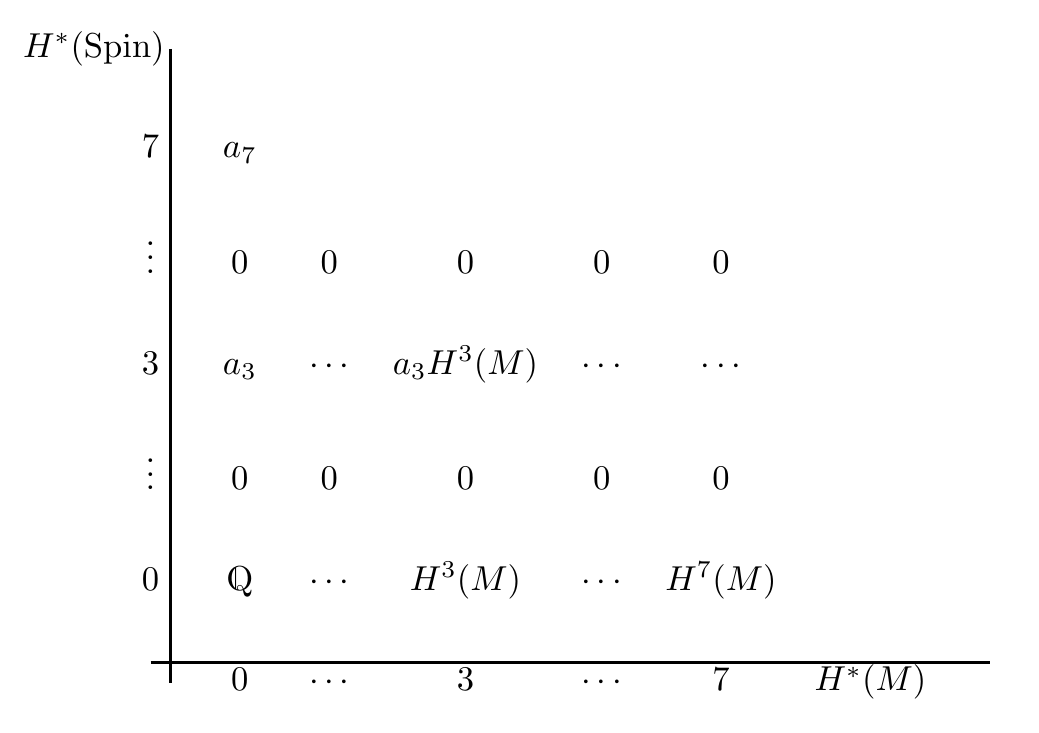}
\end{center}
As the spectral sequence converges to the cohomology of the total space and as are coefficients are $\Q$, it follows that we have a non-canonical splitting
\begin{equation*}
	H^7(Q;\Q)\cong E_\infty^{7,0}\oplus E_\infty^{4,3}\oplus E_\infty^{0,7}. 
\end{equation*}
Thus we want to calculate each of these terms.  On $E_\infty^{7,0}$, we have that the differentials $d_r$ are all zero since $E_r^{p,q}=0$ for $q<0$.  Thus the only differential of interest is $d_4:E_4^{3,3}\to E_4^{7,0}$.  Let us determine how the differential acts on generators of $E_4^{3,3}$.  Using that $E_4^{3,3}\cong E_2^{3,3}\cong \Q[a_3]\otimes H^3(M)$, then a typical generator is of the form $a_3u_3$ where $u_3\in H^3(M;\Q)$.  We also know that $d_4(a_3)=
\tfrac{1}{2} p_1$ where $p_1$ is the first Pontrjagin class of $M$, and since $M$ admits a String structure, then $p_1=0$.  Thus, since $d_r$ are derivations,
 $d_4(a_3u_3)=d(a_3)u_3+a_3d(u_3)=0$ for any generator of $E_4^{3,3}$ and 
 thus $E_\infty^{7,0}=H^7(M;\Q)$.

For $E_\infty^{0,7}$, the only relevant differentials are $d_5:E_5^{0,7}\to E_5^{5,3}$ and $d_8:E_8^{0,7}\to E_8^{8,0}$.  Since $d_r$ is zero for $r\leq4$ and $6\leq r<8$, then it follows that $E_5^{0,7}\cong \Q[a_7].$  To see what $d_5$ and $d_8$ map $\R[a_7]$ to, we will use the spectral sequence for the universal bundle $\Spin(n)\to E\Spin\to B\Spin$ along with naturality of the bundle map coming from the classifying map $f:M\to B\Spin$.  Let $F_r^{p,q}$ represent the spectral sequence for the universal bundle.  Then the map $f:M\to B\Spin$ induces maps $f^*:F_r^{p,q}\to E_r^{p,q}$ such that $f^*$ is the identity when $p=0$.  Thus $d_5(a_7)=d_5(f^*a_7)=f^*d_5(a_7)=0$ since $H^5(B\Spin;\Q)=0$ which means $F_*^{5,3}=0$.  By the same reasoning, since $d_8$ maps $a_7$ to the generator of $H^8(B\Spin;\Q)$ then for $Q$, $d_8(a_7)=\tfrac{1}{6}p_2$ where $p_2$ is the second Pontrjagin class.  Since $Q$ has a $\Five$ structure, then $d_8(a_7)=0$ and thus $E_\infty^{0,7}=E_2^{0,7}\cong \R[a_7]$.  It follows that
\begin{equation*}
	H^7(Q;\Q)\cong \Q[a_7]\oplus E_\infty^{4,3}\oplus H^7(M;\Q). 
\end{equation*}
Through a similar argument, we find that $H^7(P;\Q)\cong \Q[a_7]\oplus H^7(M;\Q)$ where we now have $E_\infty^{4,3}=0$ since $H^3(\String;\Q)=0.$  The bundle morphism $\mu:P\to Q$ induces a homomorphism $\mu^*:H^7(Q;\Q)\to H^7(P;\Q)$ and thus a homomorphism between each page of the spectral sequences.  It follows that $\mu^*$ is surjective, and that $\text{Ker}(\mu^*)=E_\infty^{4,3}$.  To finish the proof, it only remains for us to show that $E_\infty^{4,3}\cong \Q[a_3]\otimes H^4(M;\Q)$.  Indeed, the only nontrivial differential is $d_4:E_4^{4,3}\to E_4^{8,0}$ and since we have already shown that $d_4(a_3)=0$ then it follows that $d_4$ is also trivial.  
\endofproof

%
\begin{remark} {\bf (i)}
Theorem \ref{rfiv} demonstrates the degree to which the underlying Spin bundle can be used to classify lifts of the String bundles rationally.  
The difference between the integral and rational case is torsion and the Bockstein sequence corresponding to the short exact sequence
	$
	1\to \Z\to \Q\to\Q/\Z\to1
	$
prescribes to what degree that they differ.  Thus two Fivebrane structures are identified rationally if their difference corresponds to a torsion class in $H^7(M;\Z)$.

\item {\bf (ii)}
Theorem \ref{rfiv}, furthermore, gives us a similar understanding of what happens when going from Spin-Fivebrane to Fivebrane classes rationally.  Part 
{\bf (i)} tells us that every rational Spin-Fivebrane class gives rise to a rational Fivebrane class.  Part {\bf (iii)} identifies when any two rational Spin-Fivebrane classes correspond to the same rational Fivebrane class.  
\end{remark}

This says 
that rationally all the information on Fivebrane structures is essentially encoded in the underlying Spin bundles.  
In fact, it follows immediately that we have

\begin{corollary}
 If $H^4(M;\Z)$ is torsion, then the set of rational Fivebrane classes and rational 
 Spin-Fivebrane classes are in bijective correspondence.
\end{corollary} 

Settings where this occurs include the following. 

\begin{example}
{\bf (i)} The Witten manifolds $M_{k, \ell}$, 
which are $S^1$ bundles over the product of complex projective
spaces $\CC P^2 \times \CC P^1$, are classified in \cite{KS} according to 
two integers $k$ and $\ell$. They have $H^4(M_{k, \ell}; \Z)=\Z/\ell^2$. 
\item {\bf (ii)} Generalized Witten manifolds $N_{kl}$ are defined as the total spaces
of fiber bundles with fiber the lens space $L_k(\ell_2, \ell_2)$ 
and structure group $S^1$. They have $H^4(N_{kl}; \Z)\cong \Z_{|\ell_1 \ell_2|}$ \cite{E}. 
\item {\bf (iii)} Quaternionic line bundles $E$ over closed Spin manifolds of 
dimension $4k-1$ with $c_2(E) \in H^4(M; \Z)$ being torsion are considered in 
\cite{CG} via generalizations of the Kreck-Stolz invariants. 
\end{example} 

The above structures are also somewhat related to $p_i$-structures, as defined (and highlighted) 
in \cite[Def. 6.1]{9brane}.  In the absence of a rational String structure or, more 
precisely, if the obstruction $p_1^{\Q}$ for rational String structures does not vanish, 
 then the concept of a rational Spin-Fivebrane structure is equivalent to a rational 
 $p_2$-structure.   
A $p_2$-structure is a lift of $B\rm{O}$ to $B\rm{O} \langle$$p_2 \rangle$
and a rational $p_2$-structure, i.e., a $p_2^\Q$-structure,
 is a lift of the corresponding rationalizations. 
%

%

%

\subsection{Variations on rational Ninebrane classes}
\label{Sec var9}

We now extend the results from the last section to the next  higher connected cover of the 
orthogonal group $\Or$.  Following \cite{9brane}, let $2\Spin$ and $\Nine$ denote the groups 
$\Or\langle11\rangle$ and $\Or\langle15\rangle$ respectively.
Notice that in our Whitehead tower $B\Or\langle k\rangle$ for $k=10,12$ is obtained by killing homotopy groups that are completely torsion.  Hence rationally, $H^*(B\Or\langle k\rangle;\Q)\cong H^*(B\Five;\Q)$ for $k=10,12.$  So to follow along the lines of rational $\Five$ structures, we may define rational $\Nine$ structures, and so on, for all the $k$-connected covers of $\Or$ which correspond to the killing of integral homotopy groups.

\begin{definition}
A {\rm rational Ninebrane structure} is a lift of the 
$2\Spin$-principal bundle $\pi_{{}_{2\Spin}}:T\to M$ 
to the homotopy fiber $F(\tfrac{1}{240}p_3)^\Q$ of 
 the rational classifying map $f:M\to B2\Spin_\Q$.
  \end{definition}

\begin{definition}
A {\rm rational Ninebrane class} is a cohomology class $\mathcal{N}_\Q\in H^{11}(T;\Q)$
 such that $\iota_x^*\mathcal{N}=a_{11}=\tau(\tfrac{1}{240}p_3^\Q) \in H^{11}(2\Spin;\Q)$ for each inclusion $\iota_x:2\Spin\to T.$
\end{definition}

Now, just as we did in the case of Fivebrane structures, we will relate these classes to ones on the underlying Spin bundle.  In order to do this, as we compared degree 7 rational cohomology between Spin and String we need to compare the degree 11 rational cohomology of Spin and 2Spin. Letting $\rho_0$ denote the homomorphism $\rho_0:2\Spin\to\Spin,$ we consider here a principal $2\Spin$-bundle $T$ and let $Q$ again denote the induced principal $\Spin$-bundle with a bundle map $\rho:T\to Q$ which is $\rho_0$-equivariant.
\begin{lemma}\label{ison}
 The map $\rho_0:2\Spin\to \Spin$ induces an isomorphism $\rho_0^*:H^{11}(\Spin;\Q)\xrightarrow{\cong}H^{11}(2\Spin;\Q)$.
\end{lemma}
\proof
We recall from Sec. \ref{Sec rank} that the rational cohomology of Spin is given by the exterior algebra $\Exterior_\Q(x_3,x_7,x_{11},\ldots)$.  This gives a minimal model for Spin, and from the process of killing homotopy classes in the Whitehead tower, the CDGA $(\Exterior(x_{11},x_{15},\ldots),0)$ provides a minimal model for 2Spin.  Moreover, the map $\rho:2\Spin\to\Spin$ induces a map $\rho^*:(\Exterior(x_3,x_7,x_{11},\ldots),0)\to(\Exterior(x_{11},x_{15},\ldots),0)$ under which $\rho^*(x_k)=0$ for $k=3,7$ and $\rho^*(x_{11})=x_{11}$. Thus it follows that on the level of cohomology, $\rho^*:H^{11}(\Spin;\Q)\to H^{11}(2\Spin;\Q)$ is an isomorphism.
Note that, for degree reasons, $x_{11}$ generates both degree 11 cohomology groups.

\endofproof

Now  we can use Lemma \ref{ison} to relate rational Ninebrane classes to classes on the underlying Spin bundle.  
\begin{definition}
A {\rm rational Spin-Ninebrane class} is a cohomology class $\mathcal{N}_\Q$ in  $H^{11}(Q;\Q)$ such that $\iota_x^*\mathcal{N}_\Q=\tilde a_{11}\in H^{11}(\Spin;\Q)$ for each $x\in M$.
\end{definition}

We characterize these new classes as follows.

\begin{theorem}\label{rnin} Let $\pi_{2\Spin}:P\to M$ be a principal $2\Spin$-bundle with $M$ simply connected and let $\pi_\Spin:Q\to M$ be its underlying principal Spin-bundle.  
{\bf (i)} For every rational Spin-Ninebrane class $\mathcal{N}_\Q\in H^{11}(Q;\Q)$, 
the pullback $\rho^*\mathcal{N}_\Q$ is a rational Ninebrane class.
		\item {\bf (ii)} Any rational Ninebrane structure $\mathcal{M}_\Q\in H^{11}(T;\Q)$ is the image $\mathcal{M}_\Q=\rho^*\mathcal{N}_\Q$ of a rational Spin-Ninebrane class $\mathcal{N}_\Q\in H^{11}(Q;\Q)$.
		\item {\bf (iii)} Two classes $\mathcal{N}_\Q, \mathcal{N}'_\Q\in H^{11}(Q;\Q)$ will give the same rational Ninebrane structure if 
		$$
		\mathcal{N}_\Q- \mathcal{N}'_\Q=\mathcal{S}\cdot\pi_\Spin^*\psi_8+
		\mathcal{F}\cdot\pi_\Spin^*\phi_4\;,
		$$
		 where $\mathcal{S}\in H^3(Q;\Q)$ is the String structure class, $\mathcal{F}\in H^7(Q;\Q)$ is the Fivebrane structure class, $\psi_8\in H^8(M;\Q)$, and  $\phi_4\in H^4(M;\Q)$ are rational  cohomology classes.
\end{theorem}

\proof
The proof follows along similar lines as the proof of Theorem \ref{rfiv}.  Given a $2\Spin$-bundle $T$ over a manifold $M$, we have an induced $\Spin$-bundle over $M,$ by Lemma 
\ref{ison}, induced by the fibration $\rho:2\Spin\to\Spin$. 
By Lemma \ref{ison}, we also know that this fibration induces an isomorphism on rational cohomology of degree 11.  In keeping with our notation, we will denote this induced $\Spin$
bundle as $Q$.  Now as before, we will compare the Serre spectral sequences corresponding to the rational cohomology for both bundles.  As $H^k(2\Spin;\Q)=0$ for $0<k<11$, it follows easily that $H^{11}(T;\Q)=\Q[a_{11}]\oplus H^{11}(M;\Q)$.  Now for the bundle $Q$, the second page of the Serre spectral sequence is provided below.

We would like to calculate the entries $E^{p,q}_\infty$ such that $p+q=11$.  It follows immediately that $E^{2,9}_\infty=E^{3,8}_\infty=E^{5,6}_\infty=E^{6,5}_\infty=E^{7,4}_\infty=E^{9,2}_\infty=E^{10,1}_\infty=0$, and $E^{1,10}_\infty=0$ as $M$ is simply connected.  Thus
\begin{equation*}
	H^{11}(T;\Q)\cong E^{0,11}_\infty\oplus E^{4,7}_\infty
	\oplus E^{8,3}_\infty\oplus E^{11,0}_\infty.
\end{equation*}
On inspection of the universal $\Spin$ bundle, we find that $d_4(\tilde a_3)=b_4,\,d_8(\tilde a_7)=b_8,$ and $d_{12}(\tilde a_{11})=b_{12}$, where $b_i\in H^i(B\Spin;\Q)$ and $\tilde a_i\in H^i(B\Spin;\Q)$ are generators.  We also find that for all other possible differentials, $d_r(a_i)=0$.  Using functoriality of the differential maps and using the classifying map of $Q$ to compare with the universal $\Spin$ bundle, it follows that $d_r(a_3)=0$ for $r\neq4,\,d_r(a_7)=0$ for $r\neq8$, and $d_r(a_{11})=0$ for $r\neq12$.  
\begin{center}
\includegraphics[width=350px]{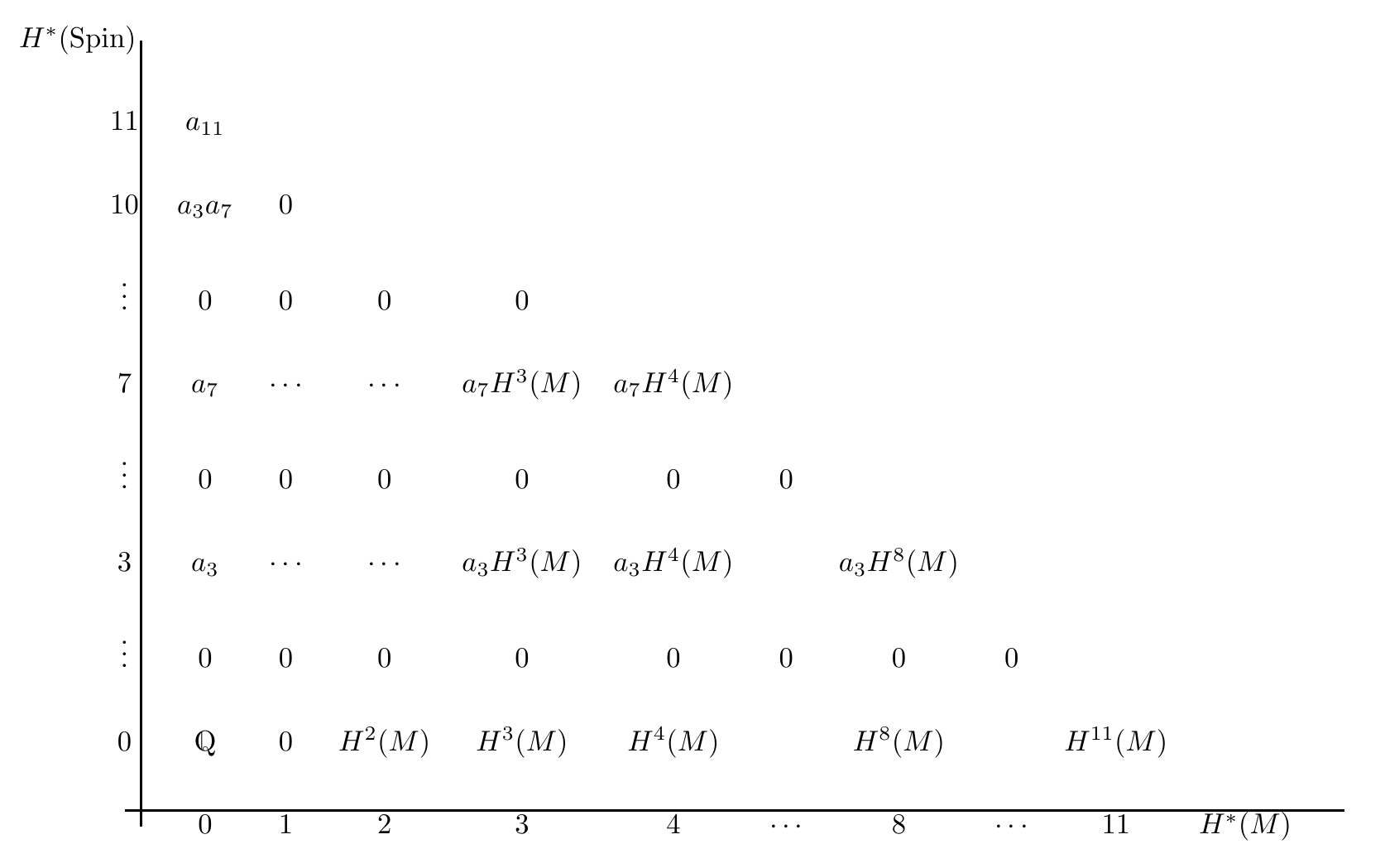}
\end{center}
Then we may proceed along the same lines as in Theorem \ref{rfiv} to identify the  following pages
\begin{align*}
&E^{0,11}_\infty\cong\Q[a_{11}]\;, \qquad \qquad
E^{4,7}_\infty\cong\Q[a_7]\otimes H^4(M;\Q)\;,\\
&E^{11,0}_\infty\cong H^{11}(M;\Q)\;, \qquad E^{8,3}_\infty\cong \Q[a_3]\otimes H^8(M;\Q)
\;,
\end{align*}
from which the theorem follows.
\endofproof

\begin{remark}
This theorem shows again that for Ninebrane structures most of the information is 
rationally encoded in the underlying Spin bundle.  While the kernel of the map which 
assigns rational Spin-Ninebrane classes to rational Ninebrane classes is larger, we 
still have a surjection.  In fact, this process should extend further to higher  structures.  
The reason is that we are making use of the fact that, rationally, there is an isomorphism
	$$
	H^{12}(B\Spin;\Q)/(p_1^\Q, p_2^\Q)\cong H^{12}(B2\Spin;\Q)\;.
	$$
Through minimal models (see Sec. \ref{Sec min}), it becomes clear that isomorphisms such as these continue to occur for higher connected covers of $\Spin$.  The difficulty in extending this definition then becomes more of a problem with determining the kernel of these maps. 
\end{remark} 

Again, we have the following.

\begin{corollary}
 If $H^4(M;\Z)$ and $H^8(M; \Z)$ are pure torsion, then the set of rational Ninebrane classes 
 and rational Spin-Ninebrane classes are in bijective correspondence.
\end{corollary}

\begin{example}
We can give a nontrivial example of a manifold $X$ which has torsion $H^4(X; \Z)$,
vanishing $H^8(X; \Z)$ and non-torsion $H^{12}(X; \Z)$. Let ${\rm SU}(2)$ be the subgroup of 
${\rm SU}(4)$ consisting of all block diagonal matrices ${\rm diag}(A, A)$ 
where $A \in {\rm SU}(2)$. Then the 12-dimensional 
quotient $X={\rm SU}(4)/{\rm SU}(2)$, viewed as 
the base of an $S^3$ bundle,  is stably parallelizable with $H^4(X; \Z)=\Z_2$, 
$H^8(X; \Z)=0$ and $H^{12}(X; \Z)=\Z$ \cite[Lemma 6.5]{Sing}.
\end{example}

%

\subsection{Gauge transformations} 
\label{Sec gauge}

We now consider automorphisms of bundles equipped with the structures that we have
just defined above. 
Let $G$ be a topological group and $G \to P \overset{\xi}{\longrightarrow} 
X$ be a continuous $G$-principal bundle. Let $\mathcal{G}(P)$ be the gauge group of 
$P$, i.e. the group of bundle automorphisms of $P$. An element 
$\eta \in \mathcal{G}(P)$ is a bundle 
isomorphism of $P$ that fits into the diagram (see e.g. \cite{Ralph})
$$
\xymatrix@R=1.5em{
P \ar[rr]^\eta_{\cong} \ar[d] && P \ar[d] 
\\
X \ar[rr]^= && X \;.
}
$$
Equivalently, $\mathcal{G}(P)$ is the group $\mathcal{P}={\rm Aut}_G(P)$ of 
$G$-equivariant homeomorphisms of $P$ covering the identity. If $P$ is the trivial bundle 
$X \times G \to X$ then $\mathcal{G}(P)$ is given by the function 
space from $X$ to $G$, i.e. 
$\mathcal{G}(P) \cong {\rm Map}(X, G)$. 
When $X$ has a basepoint $x_0\in X$, one can also consider the based gauge group 
$\mathcal{G}_0(P)$, which is the subgroup of $\mathcal{G}(P )$ whose elements fix the fiber $P_{x_0}$ ,  i.e., 
$$
\mathcal{G}_0(P)=\left\{ \eta \in \mathcal{G}(P)~|~{\rm if}~ p \in P_{x_0}
~{\rm then}~ \eta(p)=p
   \right\}\;.
$$
In relation to Fivebrane and Spin-Fivebrane classes, we consider in general a principal ${\rm O}\langle k-1\rangle$-bundle.  In the case where the structure group of this bundle lifts to ${\rm O}\langle k\rangle$ then, as noted above  in Sec. \ref{Sec var5},  these lifts are classified up to homotopy by classes in $H^{k-1}(P;\pi_{k-1}({\rm O}))$ which pull back under the fiber inclusion map $\iota_x:P_x\to P$, for every $x\in X$, to the class corresponding to a chosen generator of $H^{k-1}({\rm O}\langle k-1\rangle;\pi_{k-1}({\rm O}))$.  As gauge transformations describe homotopy equivalences 
(or even homeomorphisms) of the total space, the induced morphisms on cohomology are isomorphisms.  Gauge transformations in the based gauge group fix the fiber over the basepoint of $X$.  Thus for $\eta\in\mathcal{G}_0(P)$, we have $\eta^*\iota_{x_0}^*=\iota_{x_0}^*$.  Fixing an element in $H^*(P;\pi_{k-1}({\rm O}))$ and using that there is a canonical isomorphism between the cohomology of each fiber and the cohomology of $\Or\langle k-1\rangle$, we can view the pullback $\iota_x^*$ as an
assignment of a cohomology class in $H^*(\Or\langle k-1\rangle,\pi_{k-1}({\rm O}))$ to each element $x\in X$. 

\begin{proposition}
\label{prop-unbased}
The unbased and based gauge groups of a  ${\rm Spin}(n)\langle k \rangle_\Q$ 
bundle over $X$ are given by the mapping spaces
$$
\mathcal{G} \simeq {\rm Map}(X, \Pi_\alpha K(\Q, m_\alpha))\;, \qquad \qquad
\mathcal{G}_0 \simeq {\rm Map}_*(X, \Pi _\alpha K(\Q, m_\alpha))\;.
$$
\end{proposition} 

\proof
This follows from various classical results in the literature as well as our earlier discussion 
in Sec. \ref{Sec rat conn}. Since all of our connected cover groups are rationally abelian, 
this means that the gauge groups, which are $G$-equivariant maps,
become  simply just maps, i.e. $\mathcal{G}={\rm Map}(X, G)$. 
More precisely, the gauge transformations are $G$-equivariant homeomorphisms, 
which are equivalent to equivariant maps $P\to G$ (where $G$ acts on itself by conjugation),
and since $G$ is abelian, the map is constant on each fiber.
Alternatively, the same holds, 
by \cite[Cor 2.2]{FO}, since all components of ${\rm Map}(X, BG)$ have the same homotopy type. 
Now we use the fact that $G\simeq\Pi K(\Q, m_i)$ (see Prop. \ref{prop rat grp}).
A similar  discussion holds for the based case. 
Note that this does not require any finiteness conditions on $X$. 
\endofproof

Note that the gauge groups arise as full spaces of maps rather than 
homotopy classes of maps, in which case the gauge group would have been
some combination of cohomology classes. We will unpack some of 
these mapping spaces in order to appreciate the rich structure. We will first 
consider the more familiar ${\rm Spin}_\Q$ gauge transformations and ask
whether they lift to ${\rm String}_\Q$ gauge transformations.
To that end, consider the fibration
$K(\Q, 2) \to {\rm String}_\Q \overset{p}{\to} {\rm Spin}_\Q$
and the corresponding lift 
$$
\xymatrix{
&& {\rm String}_\Q \ar[d]^p \\
X \ar[rr]^{u_1} \ar@{..>}[rru]^{u} && {\rm Spin}_\Q\;.
}
$$
Hence we would like to consider the decomposition of  the mapping space 
${\rm Map}(X, {\rm String}_\Q)$. Given a map $u:X\to\String_\Q$, we define the mapping space 
${\rm Map}_{u}(X;{\rm String}_\Q,{\rm Spin}_\Q)$ to be the 
space of all maps $f:X\to\String_\Q$ such that $p\circ f=p\circ u=u_1$.  The fibration 
$u_1^*(p)$ is a fiber homotpically trivial fibration \cite{Thom}. Then 
${\rm Map}_{u}(X; {\rm String}_\Q, {\rm Spin}_\Q)
\simeq{\rm Map}_{u'}(X, K(\Q, 2))$ for some map 
$u': X \to K(\Q, 2)$, as ${\rm Map}_{u}(X; {\rm String}_\Q, {\rm Spin}_\Q)$
can be interpreted as a space of sections of $u_1^*(p)$, as in \cite{MoEM}. 
Then, it follows from \cite{Haef}\cite{Hansen}\cite{Thom} that the 
 function space takes the form ${\rm Map}_u(X;Y,B)=\prod_{i=0}^n K(H^{n-i}(X;G),i)$
 when $Y\to B$ is a $K(G,n)$ fibration.  Specializing to our case where we have
 $Y=\String_\Q$, $B=\Spin_\Q$, and the map $u$ corresponds to some rational $\String$ gauge transformation, 
we get the equivalence 
$$
{\rm Map}_u(X; {\rm String}_\Q, {\rm Spin}_\Q) \simeq \prod_{i=0}^2 K(H^{2-i}(X; \Q), i)\;. 
$$
This then implies the following characterization of those gauge transformations that lift.
\begin{proposition}
If $X$ is 1-connected then 
${\rm Map}_u(X; {\rm String}_\Q, {\rm Spin}_\Q)\simeq
 H^2(X; \Q) \times K(\Q, 2)$.
\end{proposition}

\begin{example}
For $S^2$ we have ${\rm Map}_u(S^2; {\rm String}_\Q, {\rm Spin}_\Q)\simeq
K(\Q, 0) \times K(\Q, 2)$,
while for $S^m$, $m >2$ we have 
${\rm Map}_u(S^m; {\rm String}_\Q, {\rm Spin}_\Q)\simeq K(\Q, 2)$. 
\end{example} 

We can similarly consider the next two fibrations
$K(\Q, 6) \to {\rm Fivebrane}_\Q \to {\rm String}_Q$
and 
 $K(\Q, 10) \to {\rm Ninebrane}_\Q \to {\rm Fivebrane}_Q$. In these cases, 
 the results in \cite{Haef}\cite{Hansen}\cite{Thom} lead us to 
\bea
{\rm Map}_u(X; {\rm Fivebrane}_\Q, {\rm String}_\Q)
& \simeq &
\prod_{i=0}^6 K(H^{6-i}(X; \Q), i)\;. 
\\
{\rm Map}_u(X; {\rm Ninebrane}_\Q, {\rm Fivebrane}_\Q) 
& \simeq &
\prod_{i=0}^{10} K(H^{10-i}(X; \Q), i)\;.
\eea 
These are considerable spaces to deal with in practice and in applications.
Nevertheless, we can get something tractable upon imposing 
some conditions.

\begin{proposition}
{\bf (i)} If $X$ is 6-connected or if $H^i(X; \Z)$ is pure torsion for $i\leq 6$, 
then given a  rational String gauge transformation $u$, the space of  lifts of $u$
to rational Fivebrane gauge transformations is given by
${\rm Map}_u(X, {\rm Fivebrane}_\Q, {\rm String}_\Q)\cong H^6(X; \Q)$.
\item {\bf (ii)} If $X$ is 10-connected or if $H^i(X; \Z)$, $i\leq 10$ is pure torsion, 
then  given a rational Fivebrane  gauge transformation $u$, the space of 
lifts of $u$ to rational Ninebrane gauge transformations is given by \newline
${\rm Map}_u(X, {\rm Ninebrane}_\Q, {\rm Fivebrane}_\Q)\cong H^{10}(X; \Q)$.
\end{proposition}

\begin{example}
The String to Fivebrane gauge transformations for the case of $S^m$ for $m \geq 7$ 
are given as $K(\Q, 6)$, while for $S^6$ they are $K(\Q, 0) \times K(\Q, 6)$. 
Similarly, the  Fivebrane to Ninebrane gauge transformations for the case of 
$S^m$ for $m \geq 11$ 
are given as $K(\Q, 10)$, while for $S^{10}$ they are $K(\Q, 0) \times K(\Q, 10)$. 
\end{example}

In terms of classifying spaces, we can consider the fibration $K(\Q,3)\to B\String_\Q\xrightarrow{Bp} B\Spin_\Q$.  Given a principal $\String_\Q$-bundle, we can consider its classifying map $f:X\to B\String_\Q$.  Then the mapping space ${\rm Map}_f(X;B\String_\Q,B\Spin_\Q)$ describes the space of all maps from $X$ to $B\String_\Q$ which lift the map $Bp\circ f:X\to B\Spin_\Q$, and we have
$$
{\rm Map}_f(X; B{\rm String_\Q}, B{\rm Spin_\Q}) \simeq \prod_{i=0}^3 K(H^{3-i}(X; \Q), i)\;. 
$$
Similarly we can again consider the fibrations related to $\Five_\Q$ and $\Nine_\Q$ 
(with the appropriate $f$) as above to obtain
\bea
{\rm Map}_f(X; B{\rm Fivebrane_\Q}, B{\rm String_\Q})
& \simeq &
\prod_{i=0}^7 K(H^{7-i}(X; \Q), i)\;. 
\\
{\rm Map}_f(X; B{\rm Ninebrane_\Q}, B{\rm Fivebrane_\Q}) 
& \simeq &
\prod_{i=0}^{11} K(H^{11-i}(X; \Q), i)\;.
\eea 
\\
In general, as the rational homotopy groups $\pi_k({\rm O})$ are $\Q$ for $k=3\text{ mod }4$, and 0 otherwise, then we can consider the set of lifts for a principal ${\rm O}\langle 4k-1\rangle_\Q$-bundle to a principal ${\rm O}\langle 4k+3\rangle_\Q$-bundle.  We have fibrations 
$K(\Q,4k-1)\to B{\rm O}\langle 4k+4\rangle_\Q \xrightarrow{\rm \xi_{4k+4}} B{\rm O}\langle 4k\rangle_\Q$. Then, using the fact that the fibrations associated with gauge transformations 
can be extended to classifying spaces \cite{FO}, which we call 
``the space of $O\langle 4k+3\rangle_{\mathbb{Q}}$ structures", we have the following. 

\begin{proposition}
\label{hoequiv}
Let $f: X \to B{\rm O}_\Q\langle 4k+4\rangle$ be a classifying map for an ${\rm O}\langle 4k+3\rangle_\Q$-bundle over $X$.  Then the space of $O\langle 4k+3\rangle_{\mathbb{Q}}$-structures on the underlying 
${\rm O}\langle 4k-1\rangle_\Q$-bundle is 
given by the space
$$
{\rm Map}_f(X; B{\rm O}\langle 4k+4\rangle_\Q, B{\rm O}\langle 4k\rangle_\Q) \simeq
 \prod_{i=0}^{4k-1} K(H^{4k-1-i}(X; \Q), i)\;. 
 $$
\end{proposition} 

Recall Proposition \ref{torsorprop} which established that the set of isomorphism classes of 
${\rm O}\langle 4k\rangle$-structures is a torsor for $H^{4k-1}(X;\pi_{4k-1}(\rm O))$.  For rational structures, 
this proposition still holds where now the set of isomorphisms classes of ${\rm O}\langle 4k+3\rangle_\Q$-structures 
is a torsor for $H^{4k-1}(X;\Q)$.  Now as isomorphic bundles have homotopic classifying maps, then we can equivalently 
interpret the set of isomorphism classes of ${\rm O}\langle 4k+3\rangle$-structures lifting an ${\rm O}\langle 4k-1\rangle$ 
bundle as the connected components 
$\pi_0$ of ${\rm Map}_f(X; B{\rm O}\langle 4k+4\rangle_\Q, B{\rm O}\langle 4k\rangle_\Q)$.  
Using the homotopy equivalence of Prop. \ref{hoequiv}, we can explicitly calculate to find
	\begin{align*}
		\pi_0\big({\rm Map}_f(X; B{\rm O}\langle 4k+4\rangle_\Q, B{\rm O}\langle 4k\rangle_\Q )\big) 
		&\cong  
		\pi_0\Big(\prod_{i=0}^{4k-1} K(H^{4k-1-i}(X; \Q), i)\Big)\\
		&\cong  H^{4k-1}(X;\Q)\;.	
	\end{align*}
This agrees with Proposition \ref{torsorprop}.

\begin{remark} {\bf (i)}
In relation to our previous discussion on rational $\Spin$-$\Five$ structures in Sec. \ref{Sec var5}, we considered 
the problem of classifying $\Five$ structures on a principal $\String$-bundle by classes on the underlying $\Spin$-bundle.  
Following along this theme, we again consider a principal $\String_\Q$-bundle and its underlying $\Spin_\Q$-bundle, 
and we consider the maps $p':\Five_\Q\to\String_\Q$ and $p:\String_\Q\to\Spin_\Q$.  Then we have the classifying 
map $f:X\to B\String_\Q$ and composition with $Bp$ gives the classifying map for the principal $\Spin_\Q$-bundle.  
Suppose further that the classifying map $f$ lifts along $Bp'$ to a map $\tilde f:X\to B\Five_\Q$.  Then the space 
${\rm Map}_{\tilde f}(X;B\Five_\Q,B\Spin_\Q)$ represents liftings of the underlying $\Spin_\Q$-bundle to 
a $\Five_\Q$ structure.  It follows from \cite{MoEM} that there is a fibration 
$$
{\rm Map}_{\tilde f}(X;B\Five_\Q,B\String_\Q)
\longrightarrow
{\rm Map}_{\tilde f}(X;B\Five_\Q,B\Spin_\Q)
\longrightarrow
{\rm Map}_{f}(X;B\String_\Q,B\Spin_\Q)
$$ 

\item {\bf (ii)} Similarly, we consider what happens with gauge transformations.  Let $u:X\to\Five_\Q$ 
be rational $\Five$ gauge transformation.  Then there is a fibration 
$$
{\rm Map}_u(X;\Five_\Q,\String_\Q)
\longrightarrow
{\rm Map}_u(X;\Five_\Q,\Spin_\Q)
\longrightarrow
{\rm Map}_{pu}(X;\String_\Q,\Spin_\Q)\;.
$$
\end{remark}

%

We now consider the rational homotopy groups of the gauge group 
$\mathcal{G}$. In particular, we will consider instances
when $\mathcal{G}$ itself admits a (variant of) rational String, Fivebrane or 
Ninebrane cover. The following results and examples can be generalized
in straightforward ways; however, we pick the dimensions indicated as these seem
to be most relevant for applications. The following few results are really corollaries 
of a slight generalization of one of the main results in \cite[Theorem 3.1]{FO}:
If $X$ has the homotopy type of a CW complex then the homotopy groups 
of the gauge group are given as
\begin{equation}\label{pigauge}
\pi_{q}^\Q(\mathcal{G}(P))\cong \sum_{r \geq 0}  H^r(X; \Q) \otimes \pi_{r+q}^\Q(G)\;.
\end{equation}
\begin{equation}
\pi_{q}^\Q(\mathcal{G}_0(P))\cong \sum_{r \geq 0}\tilde  H^r(X; \Q) \otimes \pi_{r+q}^\Q(G)\;.
\end{equation}
The statement and proof given in \cite{FO} are for $G$ a Lie group. However, 
we observe that the proof goes through for our kind of abelian topological groups;
in fact, it can also be deduced directly from Proposition \ref{prop-unbased}.
In our case, the relevant homotopy groups are $\pi_{4i-1}^\Q(G)$ for
$i \geq 1$. This then admits an interplay with $H^{4k}(X; \Q)$, somewhat 
similar to the kind of relations we encountered in the variations of the 
Spin-Fivebrane and Spin-Ninebrane structures in Sec. \ref{Sec var5} and Sec. \ref{Sec var9}. 

\begin{proposition} Let $P\to X$ be a principal $G$ bundle on an $n$-dimensional manifold $X$ 
where the group $G$ is $k$-connected.  Then the gauge group $\mathcal{G}(P)$ is $q$-connected 
where $q=k-n$.
\end{proposition}

\proof
This follows directly from the formula for the homotopy groups of the gauge group.  For the values 
$q\leq k-n$, we have $\pi^\Q_{q+r}(G)=0$ by assumption and thus $\pi^\Q_{r+q}(G)=0$ for $r\leq n$.   
\endofproof

 Specializing to the case where the group $G$ is a connected cover of ${\rm O}(n)$, we notice that as a feature 
 of Bott periodicity, the homotopy groups of the gauge group are periodic.  As we noted above, 
 Bott periodicity says that $\pi_i(\Or_\Q)=\pi_{i+4}(\Or_\Q)$ for all $i\geq 0$, and $\pi_i(\Or_\Q)=\Q$ 
 for $i=4k+3$ and 0 everywhere else.  This translates to gauge groups as follows.

\begin{proposition}\label{gaugebott}
Consider a principal $G$-bundle $P\to X$ where $G=\Or\langle k\rangle_\Q$.  The homotopy 
groups of  the gauge group $\mathcal{G}(P)$ satisfy the following periodicity conditions
	$$\pi_q(\mathcal{G}(P))=\pi_{q+4}(\mathcal{G}(P)) \qquad \qquad
{\rm and}
	\qquad \qquad
	\pi_q(\mathcal{G}_0(P))=\pi_{q+4}(\mathcal{G}_0(P))$$
for every $q\geq k$. 
\end{proposition}

\proof
For the case of $G=\Spin_\Q$, following equation \eqref{pigauge} we calculate the first four homotopy groups
	\bea
		\pi_{1}^\Q(\mathcal{G}(P))\cong&\sum_{r \geq 0}  H^{4r+2}(X; \Q) \otimes \pi_{4r+3}^\Q(G)\;,
		\qquad \qquad
		\pi_{2}^\Q(\mathcal{G}(P))\cong&\sum_{r \geq 0}  H^{4r+1}(X; \Q) \otimes \pi_{4r+3}^\Q(G)\;,\\
		\pi_{3}^\Q(\mathcal{G}(P))\cong&\sum_{r \geq 0}  H^{4r}(X; \Q) \otimes \pi_{4r+3}^\Q(G)\;, 
		\qquad \qquad \quad
		\pi_{4}^\Q(\mathcal{G}(P))\cong&\sum_{r \geq 0}  H^{4r-1}(X; \Q) \otimes \pi_{4r+3}^\Q(G)\;,\\
	\eea
	
\vspace{-3mm}	
\noindent and use Bott periodicity. 
In calculating $\pi_5$, noting that $\pi_3(G)=\pi_7(G)=\Q$, we get $\pi_5(\mathcal{G}(P))=\pi_1(\mathcal{G}(P))$.  
In general, if $L=q\text{ mod}(4)$, then $\pi_{r+q}(G)$ is non zero for $r=4i+3-L\geq 0$.  We further note that as 
$H^r(X;\Q)$ is a $\Q$-vector space then tensoring by $\Q$ induces an isomorphism.  Thus for general $q$ we have
	\begin{equation}\label{pigauge2}
		\pi_{q}^\Q(\mathcal{G}(P))\cong \sum_{r \geq 0}  
		H^{4r+3-L}(X; \Q)\otimes\pi_{4r+3-L+q}(G)\cong \sum_{r \geq 0}  H^{4r+3-L}(X; \Q),
	\end{equation}
and as $q-L=0\text{ mod(4)}$, then this equation only depends on $q$ modulo 4 and the periodicity result follows.  
Now we consider the case where $G=\Spin\langle k\rangle_\Q$.  Thus $\pi_i(G)=0$ for $i<k$.  Equation 
\eqref{pigauge2} still holds as long as $3-L+q\geq k$.  Since $3-L\geq 0$, it follows that $q\geq k$.  
The proof for the based gauge group is identical.
\endofproof
%
%

We observe that requiring the gauge group to be higher connected places conditions 
both on the underlying space $X$ as well as on lifting the structure group $G$, taken here to be a priori 
the Spin group, in the following sense. 
 
\begin{corollary} [String cover of the gauge group]
 The gauge group $\mathcal{G}$ for Spin bundles over $X$ is rationally 3-connected in the following 
cases:
\item {\bf (i)} If $\dim X \leq 3$ and $G$ is the rational String group. 
\item {\bf (ii)}  If $\dim X \leq 7$ and $G$ is the rational Fivebrane group. 
\item {\bf (iii)}  If $\dim X \leq 11$ and $G$ is the rational Ninebrane group. 
\end{corollary}

\proof
The homotopy groups through degree three of the rational gauge group from \cite{FO} 
takes the form 
$$
\sum_{i=0}^3 \pi_{i}^\Q(\mathcal{G}(P))\cong 
\bigoplus_{i=0}^3 H^i(X; \Q) \otimes \pi_3^\Q(G) 
\oplus 
\bigoplus_{i=0}^3 H^{4+i}(X; \Q) \otimes \pi_7^\Q(G) 
\oplus
\bigoplus_{i=0}^3 H^{8+i}(X; \Q) \otimes \pi_{11}^\Q(G) 
\oplus \cdots\;.
$$
The statements then follow from imposing the
condition on $\pi_i(G)$. The condition on $X$ ensures that higher terms do
not contribute. 
\endofproof

Note that
one cannot dispose of all conditions on $G$ in favor of conditions only on $X$ because
of the presence of $H^0(X; \Q)$ in the formula. 
The situation is a little improved when considering exceptional groups. 

\begin{example}
Consider an $E_8$ bundle $E$ on a manifold $X$ of dimension $n\leq12$. Since the first 
non-torsion homotopy group after $\pi_3(E_8)$ is $\pi_{15}(E_8)$ then if $X$ is 3-connected, 
the based gauge group $\mathcal{G}(E)$ is $(15-n)$-connected. It is interesting to note that, 
in particular, for M-theory extended to twelve dimensions, the based gauge group is lifted to 
its 3-connected, i.e. String, cover.
\end{example}

%
%
%
%
%

For certain nice spaces $X$ the description is more pleasant.  Consider for example the $m$-sphere $S^m$.

\begin{example} [Gauge groups of bundles over spheres] 
Consider a $G$-principal bundle $G \to P \to S^m$ over the $m$-sphere. 
Then the homotopy groups are shown in \cite{Wo} to be related as 
$\pi_n^\Q(\mathcal{G}(P)) \cong \pi_{n+m}^\Q (G) \oplus \pi_n^\Q(G)$. 
\item {\bf (i)} For $S^4$,
the gauge group  $\mathcal{G}(P)$ is rationally 3-connected if $G$ is rationally 7-connected. 


\item {\bf (ii)} For $S^8$, the gauge group $\mathcal{G}(P)$ is rationally 3-connected
 if $G$ is rationally 11-connected. 
\end{example}

%
%
%
%
%
%
We can also study the homotopy groups of the universal bundles.

\begin{example}
 Consider the universal bundle $\Spin \to E\Spin \overset{\xi_u}{\longrightarrow} B\Spin$. 
 Then, applying \cite[Thm. 4.2]{FO}, leads to 
 $\pi_k^\Q(\mathcal{G}(\xi_u))= \sum_{r\geq0} H^r(B\Spin;\Q)\otimes\pi_{r+k}(\Spin)$. 
 Now since the non-torsion generators of the cohomology of $B\Spin$ occur only in degrees that are multiples 
 of 4 and since, by transgression, the non-torsion homotopy groups are in degrees $4n+3$, it follows that 
 $\mathcal{G}(\xi_u)$ is 2-connected.  Applying Proposition \ref{gaugebott}, it follows that 
 $\pi_k^\Q(\mathcal{G}(\xi_u)=0$ if $k\neq 4n+3$.
\end{example} 



 \begin{remark}
The formulation of the rational gauge group as the rational  mapping space
 ${\rm Map}(X, G)$ comes close to the formulation via (topological or smooth) stacks. 
 Indeed, in stacks one builds a global object by starting from open subsets
 of $X$ and intersections thereof  into $G$. The group $G$ can be taken to 
 be a higher group, i.e. an $n$-group, which model the connected cover groups 
 that we consider here. However, as we discussed in the Introduction, the aim 
is to keep our tools as classical as possible in this paper. 
  \end{remark}

\vspace{.5cm}
\noindent {\large \bf Acknowledgement.} 
The authors would like to thank Corbett Redden for useful discussions and suggestions 
especially at the initial stages of this project, and Sadok Kallel, Mark Grant, Paolo Salvatore, 
and Sam Smith for useful comments. We are indebted to the anonymous 
referee for the many helpful suggestions that have significantly improved the paper.


\end{document}